\documentclass[journal]{IEEEtran}
\usepackage{amssymb}
\usepackage{epsfig}
\usepackage{amsmath}
\usepackage{graphicx}
\usepackage{subfig}

\begin{document}

\title{Stabilization for Networked Control Systems with  Simultaneous Input Delay and Markovian Packet Losses}

\author{ Hongdan Li, Chunyan Han, and Huanshui Zhang, \emph{Senior Member, IEEE}
\thanks{This work is supported by the National Natural Science
Foundation of China under Grants  61573221,
61633014, 61473134.}% <-this % stops a space
\thanks{H. Li and H. Zhang are with  School of Control Science and Engineering,
Shandong University, Jinan, Shandong, P.R.China 250061. H. Zhang is the
corresponding author.(e-mail: hszhang@sdu.edu.cn) }% <-this % stops a space
\thanks{C. Han is with School of Electrical Engineering, University of Jinan, Jinan Shandong 250022, China.}}

\markboth{}%
{Shell \MakeLowercase{\textit{et al.}}: Bare Demo of IEEEtran.cls for IEEE Journals}

\maketitle

\begin{abstract}
The mean square stabilization problem for discrete-time networked control systems (NCSs) is investigated in this article. What the difference from  most previous works is that input delay and packet losses occur simultaneously in the communication channel, moreover, the data packet dropout is modeled as a time-homogeneous Markov process which will bring some difficulties in solving the problem due to  the temporal correlation. The contributions in this paper can be summarized as two points. Firstly, the equivalence condition for the solvability of linear quadratic optimal problem in finite horizon subject to the discrete-time NCSs is expressed by solving the forward and backward stochastic difference equations (FBSDEs-M) which is derived from the maximum principle involving Markov jump and delay. Secondly, under basic assumption, the necessary and sufficient condition of mean square stabilization is given by the solutions to the coupled algebraic Riccati equations with Markov jump (CAREs-M). To our best knowledge, the problems studied in this paper are new because most previous works mainly discussed the case of only delay or packet dropout in NCSs.
\end{abstract}

% Note that keywords are not normally used for peerreview papers.

\IEEEpeerreviewmaketitle

\section{Introduction}

\IEEEPARstart{A} networked control system is defined as a control system wherein the control loops are closed through a communication network.  Actually, compared with the traditional feedback systems,  there indeed has many advantages such as reducing the weight and power, cutting the cost, improving the reliability of the system and so on, and it also has  wide applications in many fields like aircraft and high-performance automobiles, etc, see \cite{1}-\cite{5}, and references therein.  Nevertheless, in an NCS, due to the congestion in communicating channels, data packet losses and time-delay will inevitably occur, which probably cause the degradation of system performance and even instability.
In general, the study of NCSs is a meaningful but challenging subject.
In the last few years, a great deal of research on NCSs has sprung up.

Under the assumption that the packet loss is modeled as a  Bernoulli process, \cite{6} gave a necessary and sufficient condition for the stability of an NCS  where an upper bound of the packet dropout probability was given.  The output feedback control and stabilization problems for general networked control system in the case of user datagram protocol network were solved by \cite{7}. In fact, the packet dropout in the communication network is usually  modeled as either an independent identically distributed (i.i.d) Bernoulli process or a Markov chain, however,  the latter is more general and realistic.  \cite{8} mainly studied two different feedback network communication models of the update time process:  i.i.d random process and finite-state Markov chain. And the sufficient conditions of almost sure stability and mean-square stability were given for each case.  Based on a new NCS model in which sensor-to-controller and controller-to-actuator packet
dropouts history behavior were described by different independent Markov chains, \cite{9} shown  the sufficient conditions for stability by  linear matrix
inequalities  and  controller design of the NCS was given. \cite{10} presented the necessary and sufficient conditions for stability of sampled-data networked linear systems with Markovian packet losses. The stabilization for NCSs with delay has also been well studied, please see references \cite{11}-\cite{13}.

It is worth noting that most of the works focus on the case of only packet loss or delay in NCSs, there seldom concentrate on the NCSs simultaneously involved both packet loss and time-delay. In fact, as said in \cite{14}-\cite{15}, the problems that packet dropout and delay occur simultaneously are more difficult and remain challenging. However,
recently, some researches have concentrated on the NCSs  simultaneously involved both packet losses and delay, such as \cite{16} and \cite{17}. More concretely, under the assumption that packet loss and time-delay may occur simultaneously in NCSs,  \cite{16} obtained the necessary and sufficient stabilizing conditions based on the algebraic Riccati equation or Lyapunov equation. As to the NCSs with measurement packet dropout and delay, \cite{17} expressed the optimal controller with feedback gain based on a standard difference Riccati equation and an equivalent condition of stabilizing in mean square sense was given.
\begin{figure}[htbp]
  \begin{center}
  \includegraphics[width=0.3\textwidth]{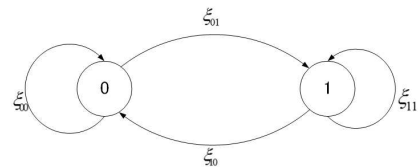}
  \caption{Two state Markovian packet dropout model} \label{fig:digit}
  \end{center}
\end{figure}

Different from \cite{16} and \cite{17} whose packet loss is molded as i.i.d Bernoulli process, in this paper we investigate the optimal LQ control and the mean square stabilization problem for discrete-time NCSs with input delay  and Markovain packet loss. It will become more complex to solve the problem due to the temporal correlation described as Fig. 1 in which transition probability is $\xi_{ij}=P(\theta_{k+1}=j|\theta_{k}=i), i,j=0,1$ and $\theta_{k}\in \{0,1\}$ denotes two state Markov chain.
Inspired by \cite{18} in which the substantial progress for the optimal LQ control has been made by solving the forward and backward difference equations (FBDEs), the necessary and sufficient conditions for the solvability of LQ optimal problem in finite horizon subject to the discrete-time NCSs is obtained by the solution to FBSDEs-M which is more difficulty to solve compared with FBDEs in \cite{18} owing to its correlations.  Further,  under exactly observability assumption,  the mean square stabilization of the NCSs can be equivalent to the positive definiteness of solutions to CAREs-M. And the main result obtained in this article can be degraded to the case of i.i.d Bernoulli packet loss, such as  \cite{16}.

 The rest of this article is mainly composed of the following sections. Section 2 gives the problem statement. Section 3 expresses the results of finite-horizon optimal control. The conclusion of stabilization is shown in section 4. In order to further illustrate the correctness of the conclusion, two numerical examples are given in section 5. The summary is provided  in section 6. There are some  relevant proofs in Appendix.

\emph{Notation }: \ ${\mathbb{R}}^n$ indicates the $n$-dimensional Euclidean space and $\mathbb{R}^{m\times n}$ denotes the norm bounded linear space of all
$m\times n$ matrices. $Y'$ is the transposition of $Y$ and if $Y\geq 0 (Y>0)$, it shows that the symmetric matrix $Y\in \mathbb{R}^{n\times n}$ is positive semi-definite(positive definite). Let a complete probability space with $\mathcal{F}_{k}$ generated by $\{\theta(0),\cdots,\theta(k)\}$ be $(\Omega,\mathcal{F}, \mathcal{F}_{k}, \mathcal{P})$. $E[\cdot|\mathcal{F}_{k}]$ means  the conditional expectation with respect
to $\mathcal{F}_{k}$ and $\mathcal{F}_{-1}$ is understood as $\{\emptyset,\Omega\}$.

\hfill %mds

\hfill %August 26, 2015
\section{Problem Statement }

Consider  the following discrete-time networked control system:
\begin{eqnarray}
x(k+1)\hspace{-3mm}&=&\hspace{-3mm}Ax(k)+\theta_{k}Bu(k-d).\label{x1}
\end{eqnarray}
The above NCS can be depicted in Fig. 2, i.e., the controller $u(k)\in {\mathbb{R}}^m$ can receive the information of the state $x(k)\in {\mathbb{R}}^n$ at time $k$. While before the control signal $u(k)$ is transmitted to plant (1), it first pass through an unreliable channel involved input time-delay $d>0$ and packet dropout which is molded as two state Markov chain $\theta_{k}\in\{0,1\}$ and its transition probability is $\xi_{ij}=\mbox{P}(\theta_{k+1}=j|\theta_{k}=i)(i,j=0,1)$.
The known initial values are $x_0, u(-1), u(-2),\cdots,u(-d)$ and the initial distribution for $\theta_{0}$ is $P(\theta_{0}=0)=q$. We assume that $\theta_{k}$ is independent of $x_0$ and $A, B$ are matrices of appropriate dimensions.
\begin{figure}[htbp]
  \begin{center}
  \includegraphics[width=0.4\textwidth]{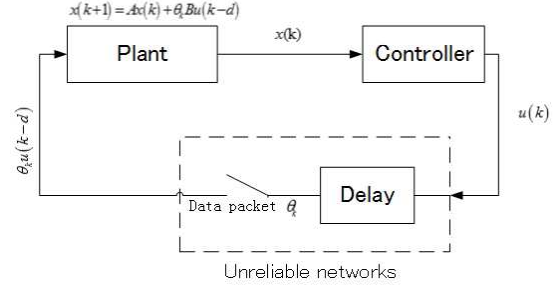}
  \caption{NCS with simultaneous input delay and packet dropout} \label{fig:digit}
  \end{center}
\end{figure}

The quadratic cost subject to  system (\ref{x1}) with infinite horizon is given by
\begin{eqnarray}
J\hspace{-3mm}&=&\hspace{-3mm}\mbox{E}\bigg\{\sum_{k=0}^{\infty}[x'(k)Qx(k)+u'(k-d)Ru(k-d)]\bigg\},\label{J2}
\end{eqnarray}
where $Q\geq 0$, $R>0$.

The following problem will be mainly discussed in this paper, i.e.,\\
\emph{ Problem 1}:  Find the $\mathcal{F}_{k-1}$-measurable controller with constant matrix gain to stabilize (\ref{x1}) while minimizing (\ref{J2}).\\
\emph{ Remark 1}:  Different from the previous works, such as \cite{16} and \cite{19}, the packet loss existing in the NCSs (\ref{x1}) is modeled as Markov process which is more general than i.i.d Bernoulli process. But due to the temporal correlation of Markov process, it will be more challenging to solve the Problem 1.

\section{Finite-horizon Optimal Control}
For discussing Problem 1, we will first introduce some associated results about the cost function with finite horizon as the following description.
\begin{eqnarray}
J_N\hspace{-3mm}&=&\hspace{-3mm}\mbox{E}\bigg\{\sum_{k=0}^Nx(k)'Qx(k)+\sum_{k=d}^Nu(k-d)'Ru(k-d)\nonumber\\
&&+x(N+1)'Hx(N+1)\bigg\},\label{J3}
\end{eqnarray}
where $N>0$ is an integer, $x(N+1)$ is the terminal state, $H$ reflects the penalty on the terminal state, the matrix functions $R\geq0$, $Q\geq0$ and $H\geq0$.

As to finite-horizon case, we will discuss  Problem 2, i.e.,

 \emph{Problem 2}: Find a $\mathcal{F}_{k-1}$-measurable controller $u(k)$ to minimize (\ref{J3}) subject to (\ref{x1}).

\emph{Lemma 1}: Problem 2 is  solvable if the following equilibrium condition is satisfied
\begin{eqnarray}
0\hspace{-3mm}&=&\hspace{-3mm}Ru(k-d)+\mbox{E}[\theta_{k}B'\lambda_k|\mathcal{F}_{k-d-1}],k=d,\cdots,N,\label{f4}
\end{eqnarray}
in which
\begin{eqnarray}
\lambda_N\hspace{-3mm}&=&\hspace{-3mm}Hx(N+1),\label{f5}\\
\lambda_{k-1}\hspace{-3mm}&=&\hspace{-3mm}Qx(k)+\mbox{E}[A'\lambda_k|\mathcal{F}_{k-1}],k=0,\cdots,N.\label{ff6}
\end{eqnarray}
\emph{Proof}: Following the results in \cite{18} and \cite{20}, the above conclusion can be similarly obtained, so we omit it.

Combining  (\ref{f4})-(\ref{ff6})  with state equation (\ref{x1}), the FBSDEs-M are established, which play a vital role in this paper.

To simplify notation, let
\begin{eqnarray*}
\prod_{\theta_{k-d}}^{\theta_{k}}\triangleq
\sum_{\theta_{k-d}=0}^{1}\xi_{\theta_{k-d-1}\theta_{k-d}}
\sum_{\theta_{k-d+1}=0}^{1}\xi_{\theta_{k-d}\theta_{k-d+1}}\cdots
\sum_{\theta_{k}=0}^{1}\xi_{\theta_{k-1}\theta_{k}},
\end{eqnarray*}
and $\xi_{\theta_{-1}0}\triangleq q, \xi_{\theta_{-1}1}\triangleq 1-q$, i.e., $\prod_{\theta_{0}}\triangleq\sum_{\theta_{0}=0}^{1}\xi_{\theta_{-1}\theta_{0}}$. And for convenience, we remark $f_{\theta_{k-1}}(k)$ as $f_{\theta_{k-1}}$.

Define the following recursions as
\begin{eqnarray}
P(N+1)\hspace{-3mm}&=&\hspace{-3mm}H,\label{f07}\\
P_{\theta_{k-1}}\hspace{-3mm}&=&\hspace{-3mm}\prod_{\theta_{k}}A'P_{\theta_{k}}A+Q-(M _{\theta_{k-1}}^{0})'
\Gamma^{-1}_{\theta_{k-1}}M_{\theta_{k-1}}^{0},\label{f7}
\end{eqnarray}
in which
\begin{eqnarray}
\Gamma_{\theta_{k-d-1}}\hspace{-3mm}&=&\hspace{-3mm}R+\prod_{\theta_{k-d}}^{\theta_{k}}\theta_{k}^{2}
B'P_{\theta_{k}}B-\sum_{s=1}^{d}\prod_{\theta_{k-d}}^{\theta_{k-s}}
[(M^{s}_{\theta_{k-s}})'\nonumber\\
\hspace{-3mm}&&\hspace{-3mm}\cdot\Gamma^{-1}_{\theta_{k-s}}M^{s}_{\theta_{k-s}}],\label{f9}\\
M^{0}_{\theta_{k-d-1}}\hspace{-3mm}&=&\hspace{-3mm}\prod_{\theta_{k-d}}^{\theta_{k}}\theta_{k}B'
P_{\theta_{k}}A^{d+1}-\sum_{s=1}^{d}\prod_{\theta_{k-d}}^{\theta_{k-s}}
[(M^{s}_{\theta_{k-s}})'\nonumber\\
\hspace{-3mm}&&\hspace{-3mm}\cdot\Gamma^{-1}_{\theta_{k-s}}M^{0}_{\theta_{k-s}}]A^{d+1-s},\label{f10}\\
M^{1}_{\theta_{k-d-1}}\hspace{-3mm}&=&\hspace{-3mm}\prod_{\theta_{k-d}}^{\theta_{k}}\theta_{k}\theta_{k-d}
B'P_{\theta_{k}}A^{d}B-\sum_{s=1}^{d}\prod_{\theta_{k-d}}^{\theta_{k-s}}\theta_{k-d}\nonumber\\
\hspace{-3mm}&&\hspace{-3mm}\cdot[(M^{s}_{\theta_{k-s}})'\Gamma^{-1}_{\theta_{k-s}}
M^{0}_{\theta_{k-s}}]A^{d-s}B,\label{f11}\\
M^{j}_{\theta_{k-d+1}}\hspace{-3mm}&=&\hspace{-3mm}\prod_{\theta_{k-d}}^{\theta_{k}}\theta_{k}\theta_{k-d+j-1}
B'P_{\theta_{k}}A^{d+1-j}B\nonumber\\
\hspace{-3mm}&&\hspace{-3mm}-\sum_{s=d-j+2}^{d}\prod_{\theta_{k-d}}^{\theta_{k-s}}[(M^{s}_{\theta_{k-s}})'
\Gamma^{-1}_{\theta_{k-s}}M^{s-d+j-1}_{\theta_{k-s}}]\nonumber\\
\hspace{-3mm}&&\hspace{-3mm}-\sum_{s=1}^{d+1-j}\prod_{\theta_{k-d}}^{\theta_{k-s}}\theta_{k-d+j-1}
[(M^{s}_{\theta_{k-s}})'\Gamma^{-1}_{\theta_{k-s}}
\nonumber\\
\hspace{-3mm}&&\hspace{-3mm}\cdot M^{0}_{\theta_{k-s}}]A^{d+1-s-j}B
, j\geq 2, \label{f12}\\
M^{s}_{\theta_{N-i-1}}\hspace{-3mm}&=&\hspace{-3mm}0, \ \ i=0, \cdots, d-1, s=0,\cdots,d.\label{f13}
\end{eqnarray}

\emph{Lemma 2}: The following relationships are established
\begin{eqnarray}
E[A'(F^{d}_{\theta_{k-1}})'|\mathcal{F}_{k-2}]\hspace{-3mm}&=&\hspace{-3mm}(M^{0}_{\theta_{k-2}})';\label{f020}\\
E[A'(F^{d-i}_{\theta_{k-i}})'|\mathcal{F}_{k-1}]\hspace{-3mm}&=&\hspace{-3mm}(F^{d-i+1}_{\theta_{k-i}})';\label{f021}\\
E[\theta_{k}B'
(F^{d-j+1}_{\theta_{k-j+1}})'|\mathcal{F}_{k-j}]\hspace{-3mm}&=&\hspace{-3mm}
(M^{j}_{\theta_{k-j}})',\label{f022}\\
\hspace{-30mm}&&\hspace{-30mm}i=1,\cdots,d-1,j=1,\cdots,d-1,\nonumber
\end{eqnarray}
where
\begin{eqnarray}
(F^{d}_{\theta_{k-1}})'\hspace{-3mm}&=&\hspace{-3mm}(S^{d}_{\theta_{k-1}})',\label{f17}\\
(F^{d-j+1}_{\theta_{k-j}})'\hspace{-3mm}&=&\hspace{-3mm}(S^{d-j+1}_{\theta_{k-1}})'
-\sum_{s=1}^{j-1}(F^{d-s+1}_{\theta_{k-s}})'\nonumber\\
\hspace{-3mm}&&\hspace{-3mm}\cdot\Gamma^{-1}_{\theta_{k-s-1}}
M^{d-j+s+1}_{\theta_{k-s-1}}, \label{f18}
\end{eqnarray}
in which
\begin{eqnarray}
(S^{1}_{\theta_{k-1}})'\hspace{-3mm}&=&\hspace{-3mm}\prod_{\theta_{k}}\theta_{k}A'P_{\theta_{k}}B
-(M^{0}_{\theta_{k-1}})'\Gamma^{-1}_{\theta_{k-1}}M^{1}_{\theta_{k-1}},\label{f19}\\
(S^{j}_{\theta_{k-1}})'\hspace{-3mm}&=&\hspace{-3mm}\prod_{\theta_{k}} A'(S^{j-1}_{\theta_{k}})'
-(M^{0}_{\theta_{k-1}})'\Gamma^{-1}_{\theta_{k-1}}M^{j}_{\theta_{k-1}}.\label{f20}
\end{eqnarray}
\emph{Proof}. See Appendix A.

Based on the above preliminaries, the solution to Problem 2 can be described as the following theorem.

\emph{Theorem 1}:
Problem 2 is uniquely solvable if and only if the difference equations (\ref{f07})-(\ref{f13}) are well defined, i.e., $\Gamma_{\theta_{k-d-1}}>0$, $k=N,\cdots,d,\ \theta_{k-d-1}=0,1$.
If this condition is satisfied, the optimal control can be given as
\begin{eqnarray}
u^{\ast}(k-d)\hspace{-3mm}&=&\hspace{-3mm}-\Gamma^{-1}_{\theta_{k-d-1}}M^{0}_{\theta_{k-d-1}}x(k-d)
-\sum_{j=1}^{d}\Gamma^{-1}_{\theta_{k-d-1}}\nonumber\\
\hspace{-3mm}&&\hspace{-3mm}\cdot M^{j}_{\theta_{k-d-1}}u(k-2d+j-1), \label{f14}
\end{eqnarray}
for $k=N,\cdots,d$.

The corresponding optimal performance index is given by
\begin{eqnarray}
J^{\ast}_N\hspace{-3mm}&=&\hspace{-3mm}\mbox{E}\bigg\{\sum_{k=0}^{d-1}x(k)'Qx(k)+x(d)'P_{\theta_{d-1}}x(d)
-x(d)'\nonumber\\
\hspace{-3mm}&&\hspace{-3mm}\cdot\sum_{s=1}^{d}(F^{d-s+1}_{\theta_{d-s}})'\Gamma^{-1}_{\theta_{d-s-1}}
E[F^{d-s+1}_{\theta_{d-s}}x(d)|\mathcal{F}_{d-s-1}]\bigg\}.\label{f15}
\end{eqnarray}
The costate $\lambda_{k-1}$ can be given as
\begin{eqnarray}
\lambda_{k-1}\hspace{-3mm}&=&\hspace{-3mm}P_{\theta_{k-1}}x(k)-\sum_{s=1}^{d}(F^{d-s+1}_{\theta_{k-s}})'
\Gamma^{-1}_{\theta_{k-s-1}}
M^{0}_{\theta_{k-s-1}}x(k-s)\nonumber\\
\hspace{-3mm}&&\hspace{-3mm}-\sum_{s=0}^{d-1}\sum_{i=d-s}^{d}
(F^{d-i+1}_{\theta_{k-i}})'\Gamma^{-1}_{\theta_{k-i-1}}
M^{s+1-d+i}_{\theta_{k-i-1}}u(k-2d+s).\nonumber\\\label{f16}
\end{eqnarray}
\emph{Proof}. See Appendix B.\\
\emph{Remark 2}: For i.i.d. Bernoulli process, i.e., $P(\theta_{k+1}=0|\theta_{k}=i)=P(\theta_{k+1}=0)=q,
P(\theta_{k+1}=1|\theta_{k}=i)=P(\theta_{k+1}=1)=1-q$,
the  recursions (\ref{f7}) can be rewritten as
\begin{eqnarray*}
P^{1}_{k}\hspace{-3mm}&=&\hspace{-3mm}A'P^{1}_{k+1}A+Q-(M^{0}_{k})'\Gamma^{-1}_{k}M^{0}_{k},\\
P^{2}_{k}\hspace{-3mm}&=&\hspace{-3mm}-M'_{k-d}\Gamma^{-1}_{k-d}M_{k-d},\\
P^{i}_{k}\hspace{-3mm}&=&\hspace{-3mm}A'P^{i-1}_{k+1}A, i=3,\cdots,d+1,
\end{eqnarray*}
where
\begin{eqnarray*}
M_{k-d}\hspace{-3mm}&=&\hspace{-3mm}(1-q)\sum^{d+1}_{j=1}B'P^{j}_{k+1}A,\\
\Gamma_{k-d}\hspace{-3mm}&=&\hspace{-3mm}R+(1-q)^{2}\sum^{d+1}_{j=1}B'P^{j}_{k+1}B+q(1-q)B'P^{1}_{k+1}B,\\
M^{0}_{k-d}\hspace{-3mm}&=&\hspace{-3mm}M_{k-d}A^{d},\\
M^{j}_{k-d}\hspace{-3mm}&=&\hspace{-3mm}(1-q)M_{k-d}A^{d-j}B, j=1,\cdots,d.
\end{eqnarray*}
Therefore, the optimal controller is as
\begin{eqnarray*}
u(k-d)&=&-\Gamma^{-1}_{k-d}M_{k-d}\{A^{d}x(k-d)+(1-q)\\
&&\cdot\sum^{d}_{j=1}A^{d-j}Bu(k-2d+j-1)\}\\
&=&-\Gamma^{-1}_{k-d}M_{k-d}E[x(k)|\mathcal{F}_{k-d-1}],
\end{eqnarray*}
which can be regarded as a special case of Theorem 1 in \cite{18}.
% needed in second column of first page if using \IEEEpubid
%\IEEEpubidadjcol

\section{Infinite-horizon Stabilizing Results}
In this section, the results of stabilization will be introduced.

\emph{Definition 1}:
The system (\ref{x1}) is mean square stabilizable if there is a ${\cal{F}}_{k-1}$-measurable controller $u(k)=T^{0}_{\theta_{k-1}}x(k)+\sum^{d}_{j=1}T^{j}_{\theta_{k-1}}u(k-d+j-1)$ satisfying $\lim_{k\rightarrow \infty}\mbox{E}[u'(k)u(k)]=0$, such that system (\ref{x1}) is asymptotically mean square stable.

To make the time horizon $N$ explicit in the finite-horizon LQR problem, we rewrite $\Gamma_{\theta_{k}}, P_{\theta_{k}}$, $P^{0}_{\theta_{k}},M^{j}_{\theta_{k}},S^{j}_{\theta_{k}},F^{j}_{\theta_{k}}$ in (\ref{f7})-(\ref{f13}) and (\ref{f17})-(\ref{f20}) as $\Gamma_{\theta_{k}}(N)$, $P_{\theta_{k}}(N)$, $P^{0}_{\theta_{k}}(N)$, $M^{j}_{\theta_{k}}(N)$, $S^{j}_{\theta_{k}}(N)$, $F^{j}_{\theta_{k}}(N)$, $j=1,\cdots,d, m=0,1$. For discussion, the terminal weight matrix $H=P_{\theta_{N}}=0$.

 \emph{Assumption 1}:
$(A,Q^{\frac{1}{2}})$ is exactly observable.

\emph{Remark 3}: The definition of exactly observable can be seen in \cite{18}.

Before the main results are discussed, we will introduce the following conclusion which will be useful to illustrate main Theorems.

\emph{Lemma 3}:  When $N\geq d$, under the condition of $R>0$,
\begin{eqnarray}
P_{\theta_{k-1}}(N)
\hspace{-1mm}-\hspace{-1mm}\sum_{s=1}^{d}[(F^{d-s+1}_{\theta_{k-s}}(N))'
\Gamma^{-1}_{\theta_{k-s-1}}(N)F^{d-s+1}_{\theta_{k-s}}(N)]
\geq 0\label{f0020}
\end{eqnarray}
  is satisfied.

\emph{Proof}. See Appendix C.

\emph{Theorem 2}:
Under Assumption 1, if the system (\ref{x1}) is mean square stabilizable, we can obtain that:

(1) \ For any $k\geq 0, m=0,1$, $P_{m}(N)$ is convergent when $N\rightarrow\infty$, i.e.,
$\lim\limits_{N\rightarrow\infty}P_{m}(N)=P_{m}$, in which $P_{m}$ satisfies
the following algebraic equations:
\begin{eqnarray}
P_{m_{d}}\hspace{-3mm}&=&\hspace{-3mm}\prod_{m_{d+1}}A'P_{m_{d+1}}A+Q-(M^{0}_{m_{d}})'\Gamma^{-1}_{m_{d}}
M^{0}_{m_{d}},\label{f22}
\end{eqnarray}
in which
\begin{eqnarray}
\Gamma_{m_{0}}\hspace{-3mm}&=&\hspace{-3mm}R+\prod_{m_{1}}^{m_{d+1}}m_{d+1}^{2}
B'P_{m_{d+1}}B\nonumber\\
\hspace{-4mm}&&\hspace{-4mm}-\sum_{s=1}^{d}\prod_{m_{1}}^{m_{s}}[(M^{d-s+1}_{m_{s}})'
\Gamma^{-1}_{m_{s}}M^{d-s+1}_{m_{s}}],\label{f24}\\
M^{0}_{m_{0}}\hspace{-3mm}&=&\hspace{-3mm}\prod_{m_{1}}^{m_{d+1}}m_{d+1}B'P_{m_{d+1}}A^{d+1}\nonumber\\
\hspace{-4mm}&&\hspace{-4mm}-\sum_{s=1}^{d}\prod_{m_{1}}^{m_{s}}[(M^{d-s+1}_{m_{s}})'
\Gamma^{-1}_{m_{s}}M^{0}_{m_{s}}]A^{s},\label{f25}\\
M^{1}_{m_{0}}\hspace{-3mm}&=&\hspace{-3mm}\prod_{m_{1}}^{m_{d+1}}m_{d+1}m_{1}B'P_{m_{d+1}}A^{d}B\nonumber\\
\hspace{-4mm}&&\hspace{-4mm}-\sum_{s=1}^{d}\prod_{m_{1}}^{m_{s}}m_{1}[(M^{d-s+1}_{m_{s}})'
\Gamma^{-1}_{m_{s}}M^{0}_{m_{s}}]A^{s-1}B,\label{f26}
\end{eqnarray}
\begin{eqnarray}
M^{j}_{m_{0}}\hspace{-3mm}&=&\hspace{-3mm}\prod_{m_{1}}^{m_{d+1}}m_{d+1}m_{j}B'P_{m_{d+1}}A^{d-j+1}B\nonumber\\
\hspace{-4mm}&&\hspace{-4mm}
-\sum_{s=j}^{d}\prod_{m_{1}}^{m_{s}}m_{j}[(M^{d-s+1}_{m_{s}})'
\Gamma^{-1}_{m_{s}}M^{0}_{m_{s}}]A^{s-j}B\nonumber\\
\hspace{-4mm}&&\hspace{-4mm}-\sum_{s=1}^{j-1}\prod_{m_{1}}^{m_{s}}
[(M^{d-s+1}_{m_{s}})'\Gamma^{-1}_{m_{s}}M^{j+1-s}_{m_{s}}].\label{f27}
\end{eqnarray}

(2) \  \begin{eqnarray*}
P_{m_{d-1}}
\hspace{-1mm}-\hspace{-1mm}\sum_{s=0}^{d-1}(F^{s+1}_{m_{s}})'
\Gamma^{-1}_{m_{s-1}}
F^{s+1}_{m_{s}}>0,
\end{eqnarray*}
in which
\begin{eqnarray}
(F^{d}_{m_{d-1}})'\hspace{-3mm}&=&\hspace{-3mm}(S^{d}_{m_{d-1}})',\label{f53}\\
(F^{d-j+1}_{m_{d-j}})'\hspace{-4mm}&=&\hspace{-4mm}(S^{d-j+1}_{m_{d-1}})'
\hspace{-2mm}-\hspace{-2mm}\sum_{s=d-j+1}^{d-1}(F^{s+1}_{m_{s}})'\Gamma^{-1}_{m_{s-1}}
M^{2d+1-j-s}_{m_{s-1}}, \label{f54}\\
(S^{1}_{m_{d-1}})'\hspace{-3mm}&=&\hspace{-3mm}\prod_{m_{d}}m_{d}A'P_{m_{d}}B
\hspace{-1mm}-\hspace{-1mm}(M^{0}_{m_{d-1}})'\Gamma^{-1}_{m_{d-1}}M^{1}_{m_{d-1}},\label{f51}\\
(S^{j}_{m_{d-1}})'\hspace{-3mm}&=&\hspace{-3mm}\prod_{m_{d}}A'(S^{j-1}_{m_{d}})'
\hspace{-1mm}-\hspace{-1mm}(M^{0}_{m_{d-1}})'\Gamma^{-1}_{m_{d-1}}M^{j}_{m_{d-1}},\label{f52}
\end{eqnarray}
$j=1,,\cdots,d$, $ m_{i}\in \{0,1\}$, $i=0,1,\cdots,d+1$.

\emph{Proof}. See Appendix D.

Next we will give the main result.

\emph{Theorem 3}: Under the condition of Assumption 1, the system (\ref{x1}) is stabilizable in the mean square sense if and only if there exists a unique solution to the Riccati-type equations (\ref{f22}) such that
\begin{eqnarray}
P_{m_{d-1}}
-\sum_{s=0}^{d-1}(F^{s+1}_{m_{s}})'
\Gamma^{-1}_{m_{s-1}}
F^{s+1}_{m_{s}}>0.\label{f0056}
\end{eqnarray}
Moreover, the optimal controller will be given as
\begin{eqnarray}
u^{\ast}(k-d)\hspace{-3mm}&=&\hspace{-3mm}-\Gamma^{-1}_{m_{0}}M^{0}_{m_{0}}x(k-d)\nonumber\\
\hspace{-3mm}&&\hspace{-3mm}-\sum_{j=1}^{d}\Gamma^{-1}_{m_{0}}M^{j}_{m_{0}}u(k-2d+j-1), \label{f56}
\end{eqnarray}
for $k\geq d$, $m_{i}\in\{0,1\}, i=0,\cdots d$.

The corresponding optimal performance index is given by
\begin{eqnarray}
J^{\ast}\hspace{-3mm}&=&\hspace{-3mm}E\bigg\{x'(0)P_{m_{1}}x(0)
-\sum_{k=0}^{d-1}u'(k-d)Ru(k-d)\nonumber\\
&&+\sum_{k=0}^{d-1}\Big[u(k-d)+
\Gamma^{-1}_{m_{d}}M^{0}_{m_{d}} x(k-d)\nonumber\\
&&+\Gamma^{-1}_{m_{d}}\sum_{j=1}^{d}M^{j}_{m_{d}}u(k-2d+j-1)\Big]'\Gamma_{m_{d}}\nonumber\\
&&\cdot\Big[u(k-d)+
\Gamma^{-1}_{m_{d}}M^{0}_{m_{d}}x(k-d)\nonumber\\
&&+\Gamma^{-1}_{m_{d}}\sum_{j=1}^{d}M^{j}_{m_{d}}u(k-2d+j-1)\Big]\bigg\}.\label{f57}
\end{eqnarray}

\emph{Proof}: See Appendix E.

\emph{Remark 4}: When $\{\theta_{k}\}$ is modeled as an i.i.d. Bernoulli process, then the special transition probability can be expressed as $P(\theta_{k+1}=0|\theta_{k}=i)=P(\theta_{k+1}=0)=q,
P(\theta_{k+1}=1|\theta_{k}=i)=P(\theta_{k+1}=1)=1-q$.
Thus, the Riccati-type equations (\ref{f22}) can be rewritten as
\begin{eqnarray*}
P^{1}&=&A'P^{1}A+Q-(M^{0})'\Gamma^{-1}M^{0},\\
P^{2}&=&-M'\Gamma^{-1}M,\\
P^{i}&=&A'P^{i-1}A, i=3,\cdots,d+1,
\end{eqnarray*}
where
\begin{eqnarray*}
M&=&(1-q)\sum^{d+1}_{j=1}B'P^{j}A,\\
M^{0}&=&MA^{d},\\
M^{j}&=&(1-q)MA^{d-j}B, j=1,\cdots,d,\\
\Gamma&=&R+(1-q)^{2}\sum^{d+1}_{j=1}B'P^{j}B+q(1-q)B'P^{1}B.
\end{eqnarray*}
Further, $F^{j}=MA^{j-1}, j=d,\cdots, 1$, in view of the relationships, (\ref{f0056}) will be reexpressed as
\begin{eqnarray*}
P^{1}+\sum^{d+1}_{j=2}P^{j}=\sum^{d+1}_{j=1}P^{j}>0.
\end{eqnarray*}
It can be seen from the results of the above transformation that the main results in this article can be degenerated to the case which  packet loss is modeled as  Bernoulli process \cite{19}.

\emph{Remark 5}:  For the delay-free case, i.e., $d=0$ in the NCSs (\ref{x1}),  the algebraic Riccati-type equations (\ref{f22}) can be reduced to the following standard  algebraic Riccati equations with Markov jump \cite{21}:
\begin{eqnarray*}
P_{m_{0}}=\prod_{m_{1}}A'P_{m_{1}}A+Q-(M^{0}_{m_{0}})'\Gamma^{-1}_{m_{0}}M^{0}_{m_{0}},
\end{eqnarray*}
where
\begin{eqnarray*}
\Gamma_{m_{0}}&=&R+\prod_{m_{1}}m^{2}_{1}B'P_{m_{1}}B,\\
M^{0}_{m_{0}}&=&\prod_{m_{1}}m_{1}B'P_{m_{1}}A,
\end{eqnarray*}
$m_{i}=0,1, i=0,1$. Moreover, (\ref{f0056}) can be expressed as $P_{m_{d-1}}>0$.

\emph{Remark 6}:  When the NCSs (\ref{x1}) exists no packet loss, i.e, $\theta_{k}=1, k=0, 1, \cdots$, the algebraic equations (\ref{f22}) can be written as:
\begin{eqnarray*}
P=A'PA+Q-(M^{0})'\Gamma^{-1}M^{0},
\end{eqnarray*}
in which
\begin{eqnarray*}
\Gamma&=&R+B'PB-\sum^{d}_{s=1}(M^{d-s+1})'\Gamma^{-1}M^{d-s+1},\\
M^{0}&=&MA^{d},\\
 M^{j}&=&MA^{d-j}B, j=1, \cdots, d,\\
M&=&B'PA-B'\sum^{d-1}_{j=0}(A^{j})'M'\Gamma^{-1}MA^{j}A.
\end{eqnarray*}
Obviously, $F^{s+1}=MA^{s},s=0,\cdots,d-1$, hence, (\ref{f0056}) can be expressed as $P-\sum^{d-1}_{s=0}(A^{s})'M'\Gamma^{-1}MA^{s}$.
Actually, the above equations are the deterministic case of (35)-(39) in \cite{18}.

\emph{Remark 7}:  Compared with \cite{22} in which optimal control problem for discrete-time MJLS with input delay was mainly investigated, in this paper, we developed the necessary and sufficient condition of the stabilization for NCSs with simultaneous input delay and Markovian packet losses.

\section{Numerical examples}
\emph{Example 1}: Considering the system (\ref{x1}) with $A=1, B=1, d=1$ with initial values $x(0)=0.1, u(-1)=-0.1$ and transition probability $\xi_{00}=0.6, \xi_{11}=0.5$ and the cost function (\ref{J2}) with $Q=R=1$. In this case, a sample path of the Markov chain $\theta_{k}$ is shown in Fig. 3 (a).
\begin{figure}[tbh]
\begin{centering}
\subfloat[A sample path with q=0.6 and p=0.5]{\begin{centering}
\includegraphics[width=4cm]{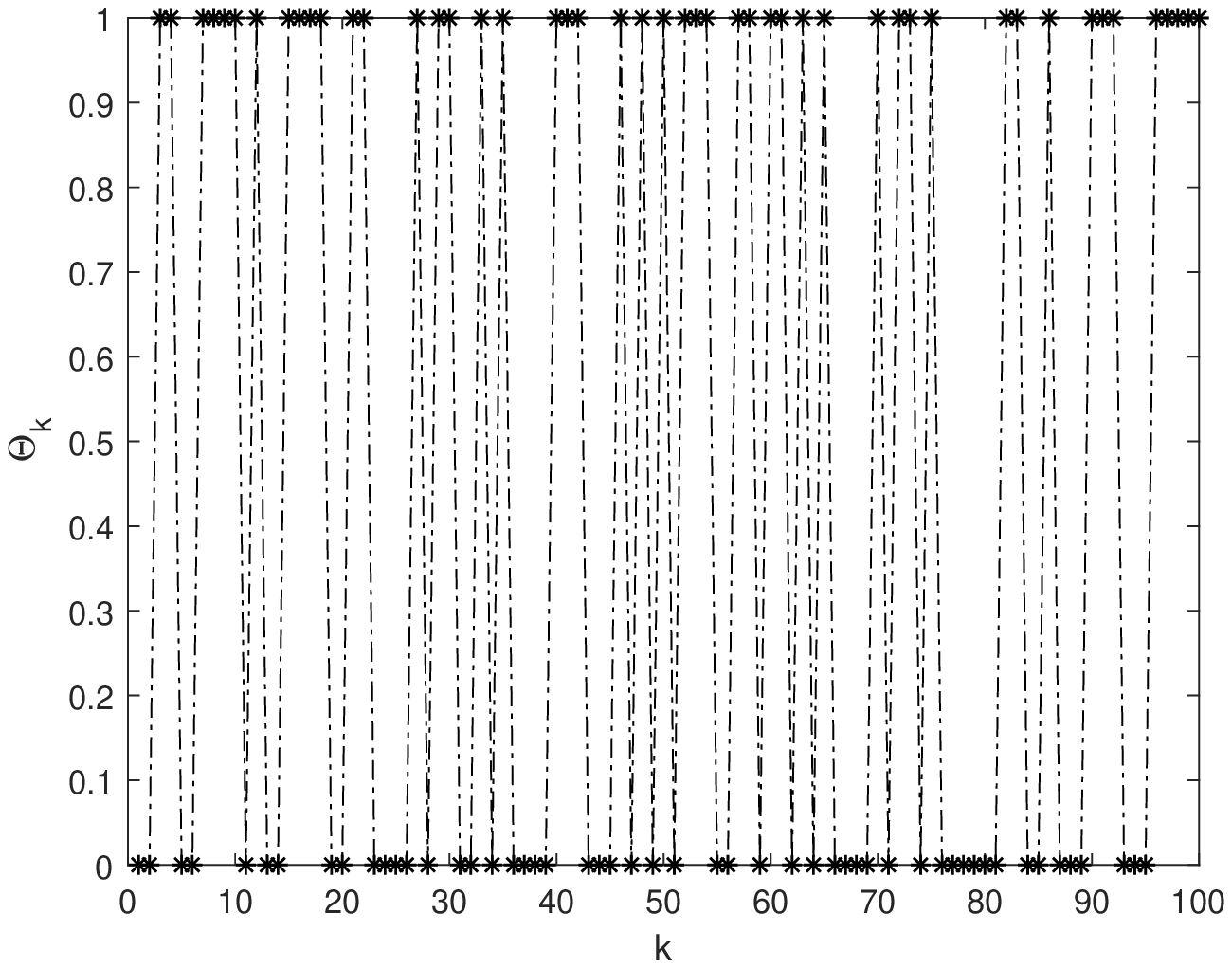}
\par\end{centering}

}\subfloat[Dynamic Behavior of $E(x_{k}'x_{k})$.]{\begin{centering}
\includegraphics[width=4cm]{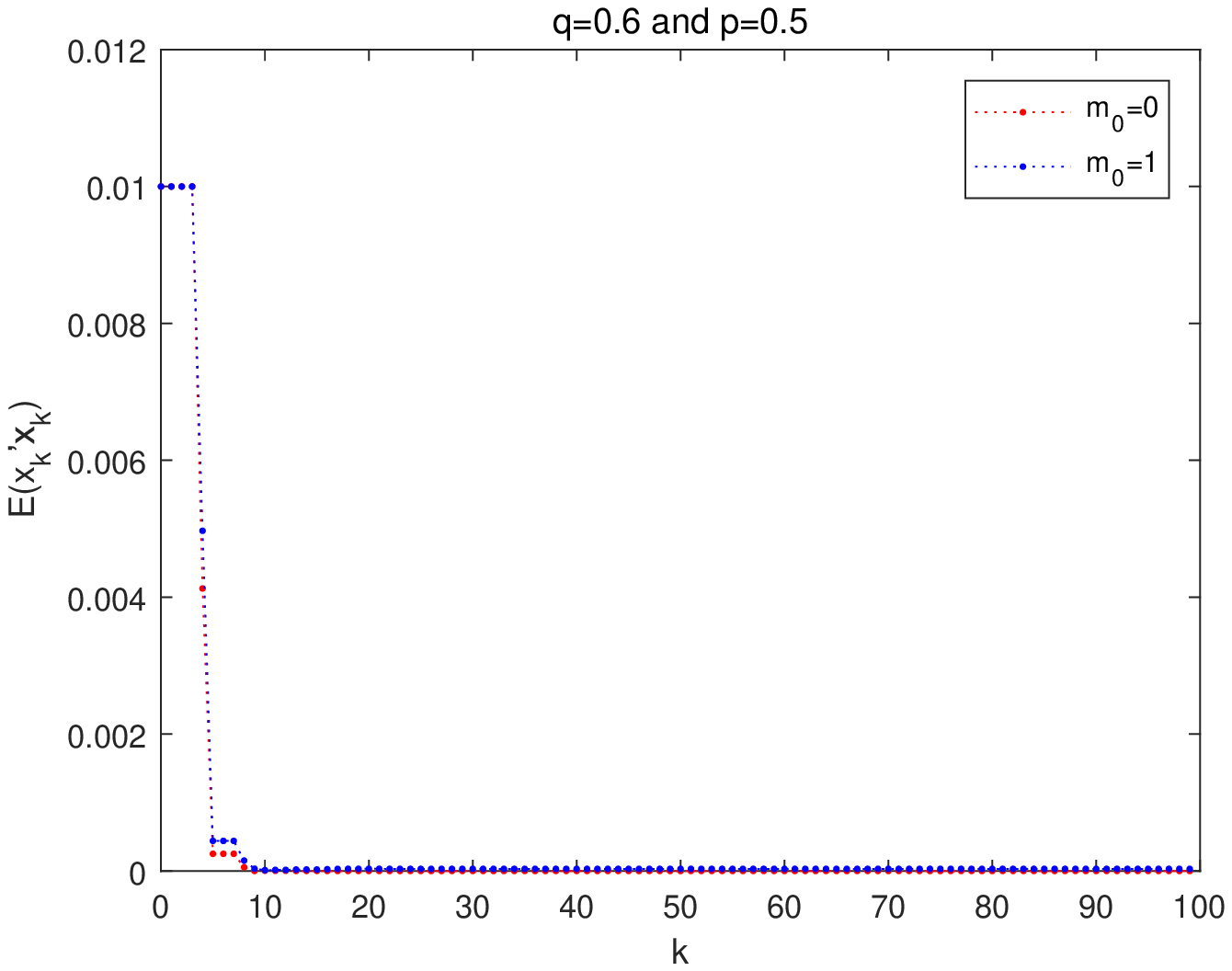}
\par\end{centering}

}
\par\end{centering}

\caption{\label{fig:PF}A sample path and $E(x_{k}'x_{k})$}
\end{figure}

It's easy to verify that Assumption 1 is satisfied. In view of Theorem 2, the following results can be obtained
$P_{0}=5.0023, P_{1}=4.9675, \Gamma_{0}=2.9538, \Gamma_{1}=2.9928, M^{0}_{0}=1.7067, M^{0}_{1}=1.7450, M^{1}_{0}=0.7746, M^{1}_{1}=0.9683, F^{1}_{0}=1.5394, F^{1}_{1}=1.9192$, further, (\ref{f0056}) can be shown as
\begin{eqnarray*}
&&P_{0}-(F^{1}_{0})'\Gamma^{-1}_{0}F^{1}_{0}=4.2000>0;\\
&&P_{0}-(F^{1}_{0})'\Gamma^{-1}_{1}F^{1}_{0}=4.2105>0;\\
&&P_{1}-(F^{1}_{1})'\Gamma^{-1}_{0}F^{1}_{1}=3.7207>0;\\
&&P_{1}-(F^{1}_{1})'\Gamma^{-1}_{1}F^{1}_{1}=3.7368>0.
\end{eqnarray*}
According to Theorem 3, the optimal controller can be expressed as
$u^{\ast}(k-1)=-0.5778x(k-1)-0.2622u(k-2)$ in the case of $m_{0}=0$ and $u^{\ast}(k-1)=-0.5908x(k-1)-0.3235u(k-2)$ in the case of $m_{0}=1, k\geq 2$. A simulation result of the designed controller is shown in Fig. 3 (b). From Fig. 3 (b), we can see that when the condition of Theorem 3 are satisfied the system is mean square stablizable.

To show the effectiveness of the result about Theorem 3, we give another example.

\emph{Example 2}: Let the coefficients in system (\ref{x1}) be taken as $A=3, B=1, d=1$ with $x(0)=0.1, u(-1)=-0.1$ and transition probability $\xi_{00}=0.9, \xi_{11}=0.7$ and the cost function (\ref{J2}) with $Q=100, R=10$. In this case, a sample path of the Markov chain $\theta_{k}$ is shown in Fig. 4 (a).
\begin{figure}[tbh]
\begin{centering}
\subfloat[A sample path with q=0.9 and p=0.7]{\begin{centering}
\includegraphics[width=4cm]{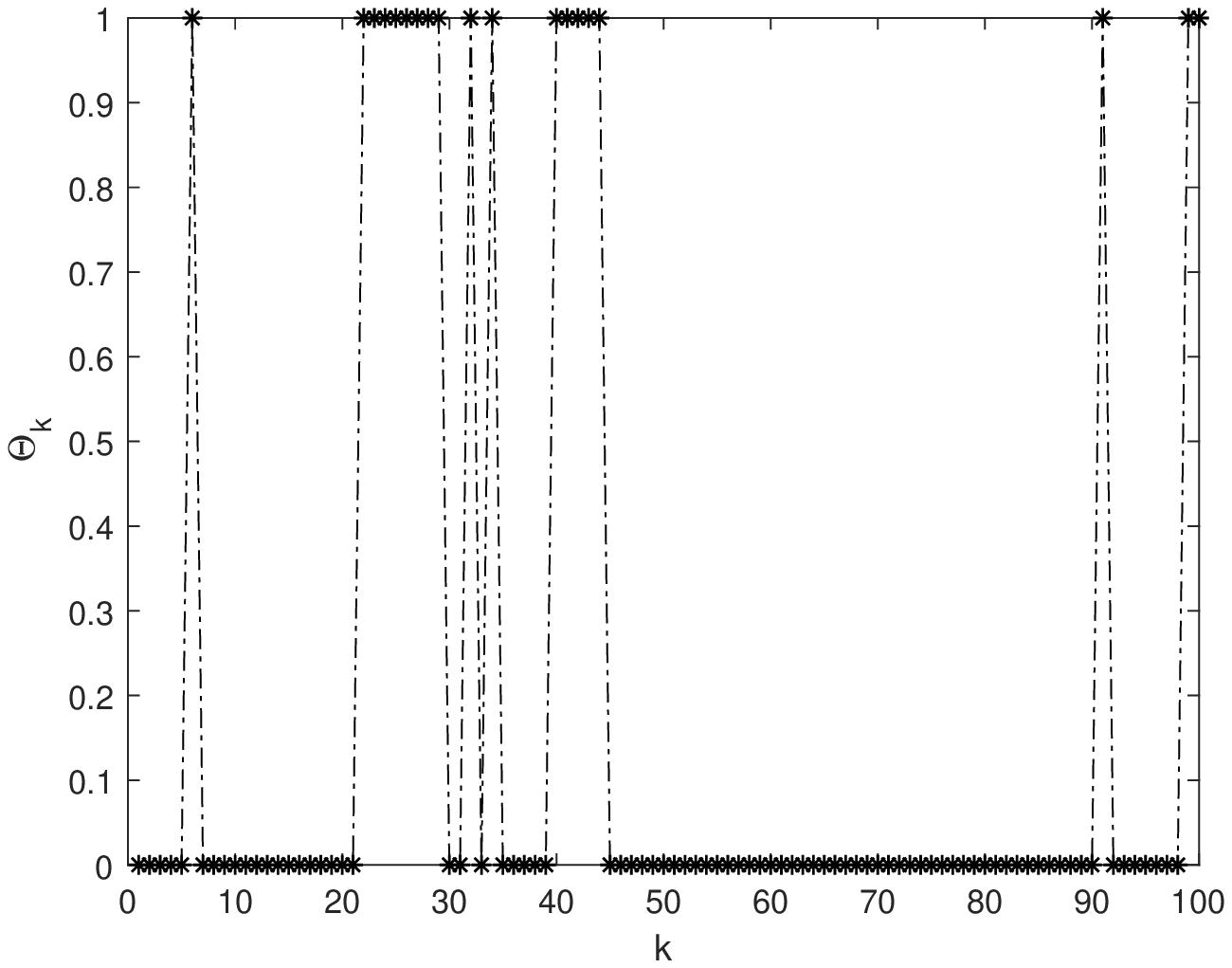}
\par\end{centering}

}\subfloat[Dynamic Behavior of $E(x_{k}'x_{k})$.]{\begin{centering}
\includegraphics[width=4cm]{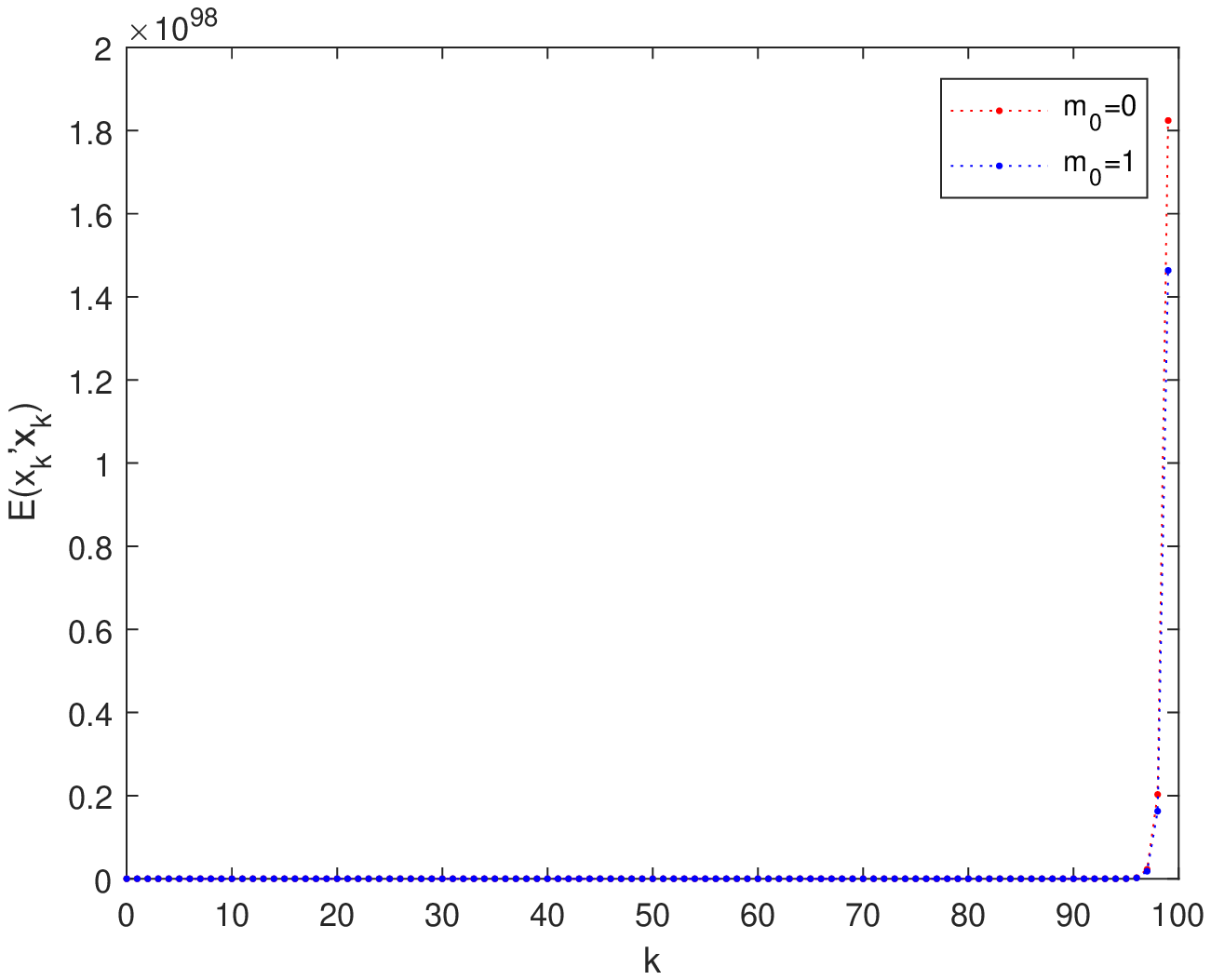}
\par\end{centering}

}
\par\end{centering}

\caption{\label{fig:PF}A sample path and $E(x_{k}'x_{k})$}
\end{figure}

   By simply calculating, the following results can be obtained
$P_{0}=5.9311, P_{1}=22.9717, \Gamma_{0}=6.8077, \Gamma_{1}=14.1673, M^{0}_{0}=-33.3076, M^{0}_{1}=-58.0541, M^{1}_{0}=-2.4371, M^{1}_{1}=-8.4956, F^{1}_{0}=-9.6267, F^{1}_{1}=-18.7327$, and
\begin{eqnarray*}
&&P_{0}-(F^{1}_{0})'\Gamma^{-1}_{0}F^{1}_{0}=-7.6819<0;\\
&&P_{0}-(F^{1}_{0})'\Gamma^{-1}_{1}F^{1}_{0}=-0.6103<0;\\
&&P_{1}-(F^{1}_{1})'\Gamma^{-1}_{0}F^{1}_{1}=-28.575<0;\\
&&P_{1}-(F^{1}_{1})'\Gamma^{-1}_{1}F^{1}_{1}=-1.7976<0.
\end{eqnarray*}
Obviously, Assumption 1 is satisfied. However, from the above values, (\ref{f0056}) is not satisfied and the corresponding controller can be obtained as
$u(k-1)=4.8926x(k-1)+0.358u(k-2)$ in the case of $m_{0}=0$ and $u^{\ast}(k-1)=4.0978x(k-1)+0.56u(k-2)$ in the case of $m_{0}=1, k\geq 2$. A simulation result of the designed controller is shown in Fig. 4 (b). It can be seen from Fig. 4 (b) that the system is not mean square stablizable.

\section{Conclusion}
In this note, we consider the optimal control and stabilization problem for discrete-time NCSs that  input delay and Markovian packet losses occur simultaneously in communication channel. Due to its temporal correlation, the analysis for such NCSs is much challenging.  The main results obtained in this paper can be  summarized as follows. Firstly, the equivalence condition for the solvability of linear quadratic optimal problem in finite horizon subject to the discrete-time NCSs is expressed by solving the FBSDEs-M. Secondly, under basic assumption, the necessary and sufficient condition of mean square stabilization is given by solutions to the CAREs-M, which can be degenerated to the case which  packet loss is modeled as  Bernoulli process, e.g., \cite{19}.

\appendices
\section{ Proof of Lemma 2}
\emph{Proof}:  Here we only prove that (\ref{f020}) is established,  (\ref{f021}) and (\ref{f022}) can be obtained similarly.

From (\ref{f17}), (\ref{f19}) and (\ref{f20})
\begin{eqnarray*}
(F^{d}_{\theta_{k-1}})'\hspace{-4mm}&=&\hspace{-4mm}(S^{d}_{\theta_{k-1}})'\\
\hspace{-4mm}&=&\hspace{-4mm}\prod\limits_{\theta_{k}}A'(S^{d-1}_{\theta_{k}})'
-(M^{0}_{\theta_{k-1}})'
\Gamma^{-1}_{\theta_{k-1}}M^{d}_{\theta_{k-1}}\\
&\vdots&\\
\hspace{-4mm}&=&\hspace{-4mm}\prod^{\theta_{k+d-2}}_{\theta_{k}}(A')^{d-1}(S^{1}_{\theta_{k+d-2}})'
-\prod^{\theta_{k+d-3}}_{\theta_{k}}(A')^{d-2}\\
\hspace{-4mm}&&\hspace{-4mm}\cdot
(M^{0}_{\theta_{k+d-1}})'\Gamma^{-1}_{\theta_{k+d-1}}M^{2}_{\theta_{k+d-1}}
-\cdots\\
\hspace{-4mm}&&\hspace{-4mm}-(M^{0}_{\theta_{k-1}})'\Gamma^{-1}_{\theta_{k-1}}
M^{d}_{\theta_{k-1}}
\end{eqnarray*}
\begin{eqnarray*}
\hspace{-4mm}&=&\hspace{-4mm}\prod^{\theta_{k+d-1}}_{\theta_{k}}\theta_{k+d-1}(A')^{d}P_{\theta_{k+d-1}}B\\
\hspace{-4mm}&&\hspace{-4mm}-\sum_{s=0}^{d-1}\Big[\prod^{\theta_{k+s-1}}_{\theta_{k}}(A')^{s}
(M^{0}_{\theta_{k+s-1}})'\Gamma^{-1}_{\theta_{k+s-1}} M^{d-s}_{\theta_{k+s-1}}\Big].
\end{eqnarray*}
On the other hand,
\begin{eqnarray*}
(M^{0}_{\theta_{k-2}})'&=&\prod^{\theta_{k+d-1}}_{\theta_{k-1}}\theta_{k+d-1}
(A')^{d+1}P_{\theta_{k+d-1}}B\\
&&-\sum_{s=1}^{d}\Big[\prod^{\theta_{k+d-1-s}}_{\theta_{k-1}}(A')^{d-s+1}
(M^{0}_{\theta_{k+d-s-1}})'\Gamma^{-1}_{\theta_{k+d-s-1}}\\
&& \cdot M^{s}_{\theta_{k+d-s-1}}\Big]\\
&=&\prod^{\theta_{k+d-1}}_{\theta_{k-1}}\theta_{k+d-1}
(A')^{d+1}P_{\theta_{k+d-1}}B\\
&&-\sum_{s=1}^{d}\Big[\prod^{\theta_{k+s-1}}_{\theta_{k-1}}(A')^{s}
(M^{0}_{\theta_{k+s-1}})'\Gamma^{-1}_{\theta_{k+s-1}} M^{d-s}_{\theta_{k+s-1}}\Big]\\
&=&E[A'(F^{d}_{\theta_{k-1}})'|\mathcal{F}_{k-2}].
\end{eqnarray*}
Hence, (\ref{f020}) is satisfied.

\section{ Proof of Theorem 1}
\emph{Proof}:  ``\emph{Necessity}" When Problem 2 has a unique solution, $\Gamma_{\theta_{k-d-1}}(k-d)>0$ will be proved by mathematical induction. To this end, define
\begin{eqnarray}
J_{k}&=&\mbox{E}\bigg\{\sum_{i=k}^N[x(i)'Qx(i)+u(i-d)'Ru(i-d)]\nonumber\\
&&+x(N+1)'Hx(N+1)\bigg\}, k=d,\cdots, N. \label{f100}
\end{eqnarray}
Considering $k=N$ in (\ref{f100}) with $x_{N}=0$, and from (\ref{x1})  we have that
\begin{eqnarray*}
J_{N}%\hspace{-3mm}&=&\hspace{-3mm}\mbox{E}[x(N)'Qx(N)+u(N-d)'Ru(N-d)+x(N+1)'Hx(N+1)]\\
\hspace{-3mm}&=&\hspace{-3mm}\mbox{E}[u(N-d)'Ru(N-d)]
+E\{E[\theta^{2}_{N}u(N-d)'B'HB\\
\hspace{-3mm}&&\hspace{-3mm}\cdot u(N-d)|\mathcal{F}_{N-1}]\}\\
\hspace{-3mm}&=&\hspace{-3mm}\mbox{E}[u(N-d)'Ru(N-d)]
+E\{E[\prod_{\theta_{N}}\theta^{2}_{N}u(N-d)'\\
\hspace{-3mm}&&\hspace{-3mm}\cdot B'HBu(N-d)|\mathcal{F}_{N-2}]\}\\
&&\vdots\\
\hspace{-3mm}&=&\hspace{-3mm}\mbox{E}[u(N-d)'Ru(N-d)]
+E\{E[\prod^{\theta_{N}}_{\theta_{N-d+1}}\theta^{2}_{N}u(N-d)'\\
\hspace{-3mm}&&\hspace{-3mm}\cdot B'HBu(N-d)|\mathcal{F}_{N-d-1}]\}\\
\hspace{-3mm}&=&\hspace{-3mm}\mbox{E}[u(N-d)'(R+\prod^{\theta_{N}}_{\theta_{N-d}}\theta^{2}_{N}B'HB)u(N-d)].
\end{eqnarray*}
Under the condition that Problem 2 has a unique solution, then for any nonzero $u(N-d)$, we have $\Gamma_{\theta_{N-d-1}}>0$.
As to $u(N-d)$, from (\ref{f4}), (\ref{f5}) and system (\ref{x1}), we obtain that
\begin{eqnarray*}
0\hspace{-3mm}&=&\hspace{-3mm}Ru(N-d)+\mbox{E}[\theta_{N}B'H(Ax(N)+\theta_{N}B\\
\hspace{-3mm}&&\hspace{-3mm}\cdot u(N-d))|\mathcal{F}_{N-d-1}]\\
\hspace{-3mm}&=&\hspace{-3mm}Ru(N-d)+\mbox{E}\{E[\theta_{N}|\mathcal{F}_{N-1}]B'HA(Ax(N-1)\\
\hspace{-3mm}&&\hspace{-3mm}+\theta_{N-1}Bu(N-d-1))+E[\theta^{2}_{N}B'HB\\
\hspace{-3mm}&&\hspace{-3mm}\cdot u(N-d)|\mathcal{F}_{N-1}]|\mathcal{F}_{N-d-1}\}
\end{eqnarray*}
\begin{eqnarray*}
\hspace{-3mm}&=&\hspace{-3mm}Ru(N-d)+\mbox{E}\{E[\prod_{\theta_{N}}\theta_{N}B'HA^{2}x(N-1)|\mathcal{F}_{N-2}]\\
\hspace{-3mm}&&\hspace{-3mm}
+E[\prod_{\theta_{N}}\theta_{N}\theta_{N-1}B'HBu(N-d-1))|\mathcal{F}_{N-2}]\\
\hspace{-3mm}&&\hspace{-3mm}+E[\prod_{\theta_{N}}\theta^{2}_{N}B'HBu(N-d)|\mathcal{F}_{N-2}]|\mathcal{F}_{N-d-1}\}\\
&&\vdots\\
\hspace{-3mm}&=&\hspace{-3mm}(R+\prod^{\theta_{N}}_{\theta_{N-d}}\theta^{2}_{N}B'HB)u(N-d)
+\prod^{\theta_{N}}_{\theta_{N-d}}\theta_{N}B'HA^{d+1}\\
\hspace{-3mm}&&\hspace{-3mm}\cdot x(N-d)
+\prod^{\theta_{N}}_{\theta_{N-d}}\theta_{N}\theta_{N-d}B'HA^{d}Bu(N-2d)\\
\hspace{-3mm}&&\hspace{-3mm}+\cdots+\prod^{\theta_{N}}_{\theta_{N-d}}\theta_{N}\theta_{N-1}B'HABu(N-d-1)\\
\hspace{-3mm}&=&\hspace{-3mm}(R+\prod^{\theta_{N}}_{\theta_{N-d}}\theta^{2}_{N}B'HB)u(N-d)
+\prod^{\theta_{N}}_{\theta_{N-d}}\theta_{N}B'HA^{d+1}\\
\hspace{-3mm}&&\hspace{-3mm}\cdot x(N-d)+\sum^{d}_{j=1}\prod^{\theta_{N}}_{\theta_{N-d}}\theta_{N}\theta_{N-j}B'HA^{j}Bu(N-d-j).
\end{eqnarray*}
Therefore, the optimal control is
\begin{eqnarray}
u(N-d)\hspace{-3mm}&=&\hspace{-3mm}-\sum_{j=1}^{d}\Gamma^{-1}_{\theta_{N-d-1}}M^{j}_{\theta_{N-d-1}}u(N-2d+j-1)\nonumber\\
\hspace{-3mm}&&\hspace{-3mm}-\Gamma^{-1}_{\theta_{N-d-1}}M^{0}_{\theta_{N-d-1}}x(N-d),\label{f101}
\end{eqnarray}
i.e., $u(N-d)$ is satisfied with (\ref{f14}) in case of $k=N$.

From (\ref{f5}), (\ref{ff6}) and (\ref{x1}), it yields that
\begin{eqnarray*}
\lambda_{N-1}\hspace{-3mm}&=&\hspace{-3mm}Qx(N)+\mbox{E}[A'H(Ax(N)+\theta_{N}Bu(N-d))|\mathcal{F}_{N-1}]\\
\hspace{-3mm}&=&\hspace{-3mm}(Q+A'HA)x(N)+\prod_{\theta_{N}}\theta_{N}A'HBu(N-d).
\end{eqnarray*}
In view of (\ref{f101}), we can see that the above formula is (\ref{f26}) with $k=N$.

Now we take any $d\leq l\leq N$, and  suppose that $\Gamma_{\theta_{k-d-1}}$ in (\ref{f9}) is positive definite, (\ref{f14}) and (\ref{f16}) are satisfied for $k\geq l+1$. Based on these assumptions,  next we will investigate that these also hold for $k=l$. Firstly, we test $\Gamma_{\theta_{l-d-1}}>0$. From (\ref{x1}), (\ref{f4}) and (\ref{ff6}), we obtain that
\begin{eqnarray*}
&&E[x(k)'\lambda_{k-1}-x(k+1)'\lambda_{k}]\\
&=&E[x(k)'Qx(k)+u(k-d)'Ru(k-d)].
\end{eqnarray*}
Adding from $k=l+1$ to $k=N$ on both sides of the above formula, when $x(l)=0$, we have that
\begin{eqnarray*}
J_{l}%\hspace{-3mm}&=&\hspace{-3mm}E[x(l+1)'\lambda_{l}]+E[x(l)'Qx(l)+u(l-d)'Ru(l-d)]\\
\hspace{-3mm}&=&\hspace{-3mm}E[\theta_{l}u(l-d)'B'\lambda_{l}]+E[u(l-d)'Ru(l-d)].
\end{eqnarray*}
%where the last equation hold with $x(l)=0$ which is similar to the discussion of $J_{N}$.

Considering (\ref{f16}), (\ref{x1}) and Lemma 1, we can obtain that
\begin{eqnarray*}
\hspace{-3mm}&&\hspace{-3mm}E[\theta_{l}u(l-d)'B'\lambda_{l}]\\
\hspace{-3mm}&=&\hspace{-3mm}E\Big\{\theta_{l}u(l-d)'B'\big[P_{\theta_{l}}(l+1)
x(l+1)-\sum_{s=1}^{d}\big((F^{d-s+1}_{\theta_{l-s+1}}(l-s+2))'\nonumber\\
\hspace{-3mm}&&\hspace{-3mm}\cdot\Gamma^{-1}_{\theta_{l-s}}(l-s+1)
E[F^{d-s+1}_{\theta_{l-s+1}}(l-s+2)x(l+1)|\mathcal{F}_{l-s}]\big)\big]\Big\}\\
\hspace{-3mm}&=&\hspace{-3mm}E\Big\{\theta_{l}u(l-d)'B'\big[P_{\theta_{l}}
\theta_{l}Bu(l-d)-\sum_{s=1}^{d}\big((F^{d-s+1}_{\theta_{l-s+1}})'\nonumber\\
\hspace{-3mm}&&\hspace{-3mm}\cdot\Gamma^{-1}_{\theta_{l-s}}
E[F^{d-s+1}_{\theta_{l-s+1}}\theta_{l}Bu(l-d)|\mathcal{F}_{l-s}]\big)\big]\Big\}\\
\hspace{-3mm}&=&\hspace{-3mm}E\Big\{u(l-d)'E[\theta^{2}_{l}B'P_{\theta_{l}}B|\mathcal{F}_{l-d-1}]u(l-d)
-u(l-d)'\\
\hspace{-3mm}&&\hspace{-3mm}\cdot\sum_{s=1}^{d}\big(E[\theta_{l}B'(F^{d-s+1}_{\theta_{l-s+1}})'
|\mathcal{F}_{l-s}]\Gamma^{-1}_{\theta_{l-s}}
E[F^{d-s+1}_{\theta_{l-s+1}}\theta_{l}B|\mathcal{F}_{l-s}]\big)u(l-d)\Big\}\\
\hspace{-3mm}&=&\hspace{-3mm}E\Big\{u(l-d)'[\prod^{l}_{l-d}\theta^{2}_{l}B'P_{\theta_{l}}B
-\sum_{s=1}^{d}\big((M^{s}_{\theta_{l+1-s}})'\Gamma^{-1}_{\theta_{l-s}}
M^{s}_{\theta_{l+1-s}}\big)]u(l-d)\Big\}.
\end{eqnarray*}
Hence, we have that
\begin{eqnarray*}
J_{l}%\hspace{-3mm}&=&\hspace{-3mm}E[\theta_{l}u(l-d)'B'\lambda_{l}]+E[u(l-d)'Ru(l-d)]\\
\hspace{-3mm}&=&\hspace{-3mm}E\Big\{u(l-d)'[R+\prod^{l}_{l-d}\theta^{2}_{l}B'P_{\theta_{l}}B
-\sum_{s=1}^{d}\big((M^{s}_{\theta_{l+1-s}})'\Gamma^{-1}_{\theta_{l-s}}
M^{s}_{\theta_{l+1-s}}\big)]\\
\hspace{-3mm}&&\hspace{-3mm}\cdot u(l-d)\Big\}\\
\hspace{-3mm}&=&\hspace{-3mm}E[u(l-d)'\Gamma_{\theta_{l-d-1}}(l-d)u(l-d)].
\end{eqnarray*}
The uniqueness of optimal control implies that $J_{l}>0$ for any nonzero $u(l-d)$. Thus, $\Gamma_{\theta_{l-d-1}}>0$.

With regard to $u(l-d)$, from (\ref{f4}), (\ref{ff6}) and Lemma 1, it yields that
\begin{eqnarray*}
0\hspace{-3mm}&=&\hspace{-3mm}Ru(l-d)+E\Big\{\theta_{l}B'[P_{\theta_{l}}x(l+1)
-\sum_{s=1}^{d}(F^{d-s+1}_{\theta_{l+1-s}})'\\
\hspace{-3mm}&&\hspace{-3mm}\cdot \Gamma^{-1}_{\theta_{l-s}}
M^{0}_{\theta_{l-s}}x(l+1-s)-\sum_{s=0}^{d-1}\sum_{i=d-s}^{d}
(F^{d-i+1}_{\theta_{l+1-i}})'\Gamma^{-1}_{\theta_{l-i}}\\
\hspace{-3mm}&&\hspace{-3mm}\cdot
M^{s+1-d+i}_{\theta_{l-i}}u(l+1-2d+s)]|\mathcal{F}_{l-d-1}\Big\}\\
\hspace{-3mm}&=&\hspace{-3mm}Ru(l-d)+E\Big\{\theta_{l}B'P_{\theta_{l}}[A^{d+1}x(l-d)
+\sum^{d}_{j=0}\theta_{l-j}\\
\hspace{-3mm}&&\hspace{-3mm}\cdot A^{j}Bu(l-d-j)]-\theta_{l}B'[(F^{d}_{\theta_{l}})'
\Gamma^{-1}_{\theta_{l-1}}
M^{0}_{\theta_{l-1}}x(l)\\
\hspace{-3mm}&&\hspace{-3mm}+(F^{d-1}_{\theta_{l-1}})'
\Gamma^{-1}_{\theta_{l-2}}
M^{0}_{\theta_{l-2}}x(l-1)+\cdots+(F^{1}_{\theta_{l+1-d}})'
\Gamma^{-1}_{\theta_{l-d}}\\
\hspace{-3mm}&&\hspace{-3mm}\cdot
M^{0}_{\theta_{l-d}}x(l+1-d)]-\sum_{s=0}^{d-1}\sum_{i=d-s}^{d}(M^{i}_{\theta_{l-i}})'\Gamma^{-1}_{\theta_{l-i}}
M^{s+1-d+i}_{\theta_{l-i}}\\
\hspace{-3mm}&&\hspace{-3mm}\cdot u(l+1-2d+s)|\mathcal{F}_{l-d-1}\Big\}
\end{eqnarray*}
\begin{eqnarray*}
\hspace{-3mm}&=&\hspace{-3mm}Ru(l-d)+E\Big\{\theta_{l}B'P_{\theta_{l}}[A^{d+1}x(l-d)
+\sum^{d}_{j=0}\theta_{l-j}\\
\hspace{-3mm}&&\hspace{-3mm}\cdot A^{j}Bu(l-d-j)]-[(M^{1}_{\theta_{l-1}})'
\Gamma^{-1}_{\theta_{l-1}}
M^{0}_{\theta_{l-1}}(A^{d}x(l-d)\\
\hspace{-3mm}&&\hspace{-3mm}
+\sum^{d}_{j=1}\theta_{l-j}A^{j-1}Bu(l-d-j))+\cdots+(M^{d}_{\theta_{l-d}})'
\Gamma^{-1}_{\theta_{l-d}}\\
\hspace{-3mm}&&\hspace{-3mm}\cdot
M^{0}_{\theta_{l-d}}(Ax(l-d)
+\theta_{l-d}Bu(l-2d))]-\sum_{s=0}^{d-1}\sum_{i=d-s}^{d}(M^{i}_{\theta_{l-i}})'\nonumber\\
\hspace{-3mm}&&\hspace{-3mm}\cdot\Gamma^{-1}_{\theta_{l-i}}
M^{s+1-d+i}_{\theta_{l-i}}u(l+1-2d+s)|\mathcal{F}_{l-d-1}\Big\}\\
\hspace{-3mm}&=&\hspace{-3mm}\Gamma_{\theta_{l-d-1}}u(l-d)+\Big[\prod_{\theta_{l-d}}^{\theta_{l}}\theta_{l}B'
P_{\theta_{l}}A^{d+1}-\sum_{s=1}^{d-1}\prod_{\theta_{l-d}}^{\theta_{l-s}}
\Big((M^{s}_{\theta_{l-s}})'\\
\hspace{-3mm}&&\hspace{-3mm}\cdot \Gamma^{-1}_{\theta_{l-s}}M^{0}_{\theta_{l-s}}\Big)
A^{d+1-s}\Big]x(l-d)+\sum_{j=1}^{d}\Big[\prod_{\theta_{l-d}}^{\theta_{l}}\theta_{l}\theta_{l-d+j-1}\\
\hspace{-3mm}&&\hspace{-3mm}\cdot
B'P_{\theta_{l}}A^{d+1-j}B-\sum_{s=1}^{d+1-j}\prod_{\theta_{l-d}}^{\theta_{l-s}}
\theta_{l-d+j-1}[(M^{s}_{\theta_{l-s}})'\Gamma^{-1}_{\theta_{l-s}}\\
\hspace{-3mm}&&\hspace{-3mm}\cdot
M^{0}_{\theta_{l-s}}]A^{d+1-s-j}B-\sum_{s=d-j+2}^{d}\prod_{\theta_{l-d}}^{\theta_{l-s}}[(M^{s}_{\theta_{l-s}})'\\
\hspace{-3mm}&&\hspace{-3mm}\cdot
\Gamma^{-1}_{\theta_{l-s}}lM^{s-d+j-1}_{\theta_{l-s}}]\Big]u(l-2d+j-1).
\end{eqnarray*}
Therefore,
\begin{eqnarray}
u(l-d)\hspace{-3mm}&=&\hspace{-3mm}-\sum_{j=1}^{d}\Gamma^{-1}_{\theta_{l-d-1}}M^{j}_{\theta_{l-d-1}}u(l-2d+j-1)\nonumber\\
\hspace{-3mm}&&\hspace{-3mm}-\Gamma^{-1}_{\theta_{l-d-1}}M^{0}_{\theta_{l-d-1}}x(l-d),\label{f102}
\end{eqnarray}
i.e., (\ref{f14}) is established with $k=l$.

As to the expression of $\lambda_{l-1}$, from (\ref{x1}), (\ref{ff6}) and (\ref{f102}), we can obtain that (\ref{f16}) hold with $k=l$ in a similar calculating way with $u(l-d)$.

``\emph{Sufficiency}" Under the condition  that $\Gamma_{\theta_{k-d-1}}>0$, the unique solvability of Problem 2 will be proved. To this end, let
\begin{eqnarray}
\mathcal{L}(k)\hspace{-3mm}&=&\hspace{-3mm}E\bigg\{x'(k)P_{\theta_{k-1}}x(k)
-x'(k)\sum_{s=1}^{d}(F^{d-s+1}_{\theta_{k-s}})'
\Gamma^{-1}_{\theta_{k-s-1}}\nonumber\\
\hspace{-3mm}&&\hspace{-3mm}\cdot
M^{0}_{\theta_{k-s-1}}x(k-s)-x'(k)\sum_{s=1}^{d-1}\sum_{i=d-s}^{d}
(F^{d-i+1}_{\theta_{k-i}})'\Gamma^{-1}_{\theta_{k-i-1}}\nonumber\\
\hspace{-3mm}&&\hspace{-3mm}\cdot
M^{s+1-d+i}_{\theta_{k-i-1}}u(k-2d+s)\bigg\}.\label{f104}
\end{eqnarray}
From Lemma 2 and (\ref{x1}), we can obtain that
\begin{eqnarray}
\hspace{-3mm}&&\hspace{-3mm}\mathcal{L}(k)-\mathcal{L}(k+1)\nonumber\\
\hspace{-3mm}&=&\hspace{-3mm}E\Big\{x(k)'[P_{\theta_{k-1}}
-\prod_{\theta_{k}}A'P_{\theta_{k}}A+E[A'(F^{d}_{\theta_{k-1}})'|\mathcal{F}_{k-1}]
\nonumber\\
\hspace{-3mm}&&\hspace{-3mm}\cdot \Gamma^{-1}_{\theta_{k-1}}M_{\theta_{k-1}}^{0}]x(k)+u(k-d)'\big[-\prod_{\theta_{k-d}}^{\theta_{k}}\theta_{k}^{2}
B'P_{\theta_{k}}B\nonumber\\
\hspace{-3mm}&&\hspace{-3mm}-\sum_{s=1}^{d}\prod_{\theta_{k-d}}^{\theta_{k-s}}
[E[\theta_{k}B'F^{d+1-s}_{\theta_{k-s+1}}|\mathcal{F}_{k-s}]
\Gamma^{-1}_{\theta_{k-s}}M^{s}_{\theta_{k-s}}]\big]u(k-d)\nonumber
\end{eqnarray}
\begin{eqnarray}
\hspace{-3mm}&&\hspace{-3mm}
+u(k-d)'M^{0}_{\theta_{k-d-1}}x(k-d)+u(k-d)'\sum_{j=1}^{d}M^{j}_{\theta_{k-d-1}}\nonumber\\
\hspace{-3mm}&&\hspace{-3mm}\cdot u(k-2d+j-1)
+x(k-d)'(M^{0}_{\theta_{k-d-1}})'u(k-d)\nonumber\\
\hspace{-3mm}&&\hspace{-3mm}+x(k-d)'(M^{0}_{\theta_{k-d-1}})'\Gamma^{-1}_{\theta_{k-d-1}}
M^{0}_{\theta_{k-d-1}}x(k-d)
+x(k-d)'\nonumber\\
\hspace{-3mm}&&\hspace{-3mm}\cdot(M^{0}_{\theta_{k-d-1}})'\Gamma^{-1}_{\theta_{k-d-1}}
\sum_{j=1}^{d}M^{j}_{\theta_{k-d-1}}u(k-2d+j-1)\nonumber\\
\hspace{-3mm}&&\hspace{-3mm}+\sum_{j=1}^{d}u(k-2d+j-1)'(M^{j}_{\theta_{k-d-1}})'u(k-d)+\sum_{j=1}^{d}u(k-2d+j-1)'\nonumber\\
\hspace{-3mm}&&\hspace{-3mm}\cdot(M^{j}_{\theta_{k-d-1}})'\Gamma^{-1}_{\theta_{k-d-1}}
M^{0}_{\theta_{k-d-1}}x(k-d)+\big(\sum_{j=1}^{d}u(k-2d+j-1)\nonumber\\
\hspace{-3mm}&&\hspace{-3mm}\cdot M^{j}_{\theta_{k-d-1}}\big)'
\Gamma^{-1}_{\theta_{k-d-1}}\big(\sum_{j=1}^{d}u(k-2d+j-1)
M^{j}_{\theta_{k-d-1}}\big)\Big\}\nonumber\\
\hspace{-3mm}&=&\hspace{-3mm}E\Big\{x(k)'Qx(k)+u(k-d)'Ru(k-d)-\Big(u(k-d)
+\Gamma^{-1}_{\theta_{k-d-1}}\nonumber\\
\hspace{-3mm}&&\hspace{-3mm}\cdot M^{0}_{\theta_{k-d-1}}x(k-d)
+\sum_{j=1}^{d}\Gamma^{-1}_{\theta_{k-d-1}}M^{j}_{\theta_{k-d-1}}u(k-2d+j-1)\Big)'\nonumber\\
\hspace{-3mm}&&\hspace{-3mm}\cdot
\Gamma_{\theta_{k-d-1}}\Big(u(k-d)+\Gamma^{-1}_{\theta_{k-d-1}}M^{0}_{\theta_{k-d-1}}x(k-d)\nonumber\\
\hspace{-3mm}&&\hspace{-3mm}+\sum_{j=1}^{d}\Gamma^{-1}_{\theta_{k-d-1}}M^{j}_{\theta_{k-d-1}}u(k-2d+j-1)\Big)\Big\}.\label{f1105}
\end{eqnarray}
Summing up from $k=d$ to $k=N$ on both sides of (\ref{f1105}), and in view of $\Gamma_{\theta_{k-d-1}}>0$ for $k\geq d$, therefore, the optimal controller and optimal cost can be given by (\ref{f14}) and (\ref{f15}). The sufficiency proof is completed.

%\begin{eqnarray*}
%J_{N}\hspace{-3mm}&=&\hspace{-3mm}\mathcal{L}(d)+E\Big\{\sum_{k=0}^{d-1}x(k)'Qx(k)+\sum_{k=d}^{N}\Big(u(k-d)\\
%\hspace{-3mm}&&\hspace{-3mm}+\Gamma^{-1}_{\theta_{k-d-1}}M^{0}_{\theta_{k-d-1}}x(k-d)\\
%\hspace{-3mm}&&\hspace{-3mm}+\sum_{j=1}^{d}\Gamma^{-1}_{\theta_{k-d-1}}M^{j}_{\theta_{k-d-1}}u(k-2d+j-1)\Big)'
%\Gamma_{\theta_{k-d-1}}\nonumber\\
%\hspace{-3mm}&&\hspace{-3mm}\cdot\Big(u(k-d)+\Gamma^{-1}_{\theta_{k-d-1}}M^{0}_{\theta_{k-d-1}}x(k-d)\nonumber\\
%\hspace{-3mm}&&\hspace{-3mm}+\sum_{j=1}^{d}\Gamma^{-1}_{\theta_{k-d-1}}M^{j}_{\theta_{k-d-1}}u(k-2d+j-1)\Big)\Big\}.\label{f105}
%\end{eqnarray*}
%In view of $\Gamma_{\theta_{k-d-1}}(k-d)>0$ for $k\geq d$, therefore, the optimal controller and optimal cost can be given by (\ref{f14}) and (\ref{f15}). The sufficiency proof is completed.

\section{Proof of Lemma 3}
\emph{Proof}: Let cost (\ref{J3}) subject to system (\ref{x1}) start at $d$ and end at $m, m\geq N$ and denote it as $\hat{J}_{d}(m)$. Following from Lemma 1 in \cite{18}, we derive that when $R>0$, Problem 2 has a unique solution. Thus, recall the conclusion of Theorem 1, the optimal value can be expressed as
\begin{eqnarray}
\hat{J}^{\ast}_{d}(m)\hspace{-3mm}&=&\hspace{-3mm}\mbox{E}\bigg\{x(d)'P_{\theta_{d-1}}
x(d)
-x(d)'\sum_{s=1}^{d}(F^{d-s+1}_{\theta_{d-s}})'\nonumber\\
\hspace{-3mm}&&\hspace{-3mm}\cdot\Gamma^{-1}_{\theta_{d-s-1}}
E[F^{d-s+1}_{\theta_{d-s}}x(d)|\mathcal{F}_{d-s-1}]\bigg\}\nonumber\\
\hspace{-3mm}&=&\hspace{-3mm}x(d)'\Big\{P_{\theta_{d-1}}-\sum_{s=1}^{d}
(F^{d-s+1}_{\theta_{d-s}})'\nonumber\\
\hspace{-3mm}&&\hspace{-3mm}\cdot\Gamma^{-1}_{\theta_{d-s-1}}
F^{d-s+1}_{\theta_{d-s}}\Big\}x(d) \geq 0. \label{f0023}
\end{eqnarray}

The arbitrary of $x(d)$ yields that
\begin{eqnarray}
\hspace{-1mm}&&\hspace{-5mm}P_{\theta_{d-1}}(m)\hspace{-1mm}-\hspace{-1mm}\sum_{s=1}^{d}(F^{d-s+1}_{\theta_{d-s}}(m))'
\Gamma^{-1}_{\theta_{d-s-1}}(m)
F^{d-s+1}_{\theta_{d-s}}(m)\nonumber\\
&\geq&0.\label{f0021}
\end{eqnarray}
Let $m=N-k+d,k\geq d$. Noting the time-variance and (\ref{f0021}), we know that
\begin{eqnarray}
\hspace{-1mm}&&\hspace{-5mm}P_{\theta_{k-1}}(N)
\hspace{-1mm}-\hspace{-1mm}\sum_{s=1}^{d}[(F^{d-s+1}_{\theta_{k-s}}(N))'
\Gamma^{-1}_{\theta_{k-s-1}}(N)F^{d-s+1}_{\theta_{k-s}}(N)]\nonumber\\
&\geq&0.\label{f0022}
\end{eqnarray}

\section{ Proof of Theorem 2}

\emph{Proof}: Firstly, we will illustrate the convergence of $P_{m}(k,N), m=0,1$.  To this end, denote
\begin{eqnarray*}
z(k)\hspace{-3mm}&=&\hspace{-3mm}\left[
  \begin{array}{c}
x(k)\\
u(k-1)\\
\vdots\\
u(k-d)
  \end{array}
\right],\tilde{A}_{k}\hspace{-1mm}=\hspace{-1mm}\left[
  \begin{array}{ccccc}
A&0&\cdots&0&\theta_{k}B\\
0&0&\cdots&0&0\\
0&I&\cdots&0&0\\
\vdots&\vdots&\ddots&\vdots&\vdots\\
0&0&\cdots&I&0
  \end{array}
\right],\\
\tilde{B}'\hspace{-3mm}&=&\hspace{-3mm}\left[
  \begin{array}{ccccc}
0&I&0&\cdots&0
\end{array}
\right],
\end{eqnarray*}
hence, system (\ref{x1}) can be expressed as the following Markov jump linear system (MJLS)
\begin{eqnarray}
z(k+1)=\tilde{A}_{k}z(k)+\tilde{B}u(k).\label{f36}
\end{eqnarray}
 Considering the infinite cost function subject to (\ref{f36}) as follows:
\begin{eqnarray}
\tilde{J}=\sum_{k=0}^{\infty}E[z'(k)\tilde{Q}z(k)+u'(k)Ru(k)],\label{f38}
\end{eqnarray}
in which $\tilde{Q}=\left[
  \begin{array}{ccccc}
Q&  & &\\
 & 0 & &\\
 &  & \ddots & \\
 &  &  & 0
 \end{array}
\right].$
The corresponding finite cost is
\begin{eqnarray}
\tilde{J}(N)\hspace{-3mm}&=&\hspace{-3mm}\sum_{k=0}^{N}E[z'(k)\tilde{Q}z(k)+u'(k)Ru(k)\nonumber\\
\hspace{-3mm}&&\hspace{-3mm}+z'(N+1)\tilde{H}z(N+1)].\label{f39}
\end{eqnarray}
By Maximum principle, the associated costate equation and equilibrium condition can be expressed as follows
\begin{eqnarray}
\left\{
\begin{array}{lll}
\beta_{k-1}=\tilde{Q}z(k)+E[(\tilde{A}+\theta_{k}\tilde{B}_{0})'\beta(k)|\mathcal{F}_{k-1}],\\
\beta_{N}=\tilde{H}z(N+1),\\
0=Ru(k)+E[\tilde{B}'\beta_{k}|\mathcal{F}_{k}],
\end{array}
\right.\label{f40}
\end{eqnarray}
in which
\begin{eqnarray*}
\tilde{A}_{k}&=&\left[
  \begin{array}{ccccc}
A&0&\cdots&0&0\\
0&0&\cdots&0&0\\
0&I&\cdots&0&0\\
\vdots&\vdots&\ddots&\vdots&\vdots\\
0&0&\cdots&I&0
  \end{array}
\right]+\theta_{k}\left[
  \begin{array}{ccccc}
0&\cdots&0&B\\
0&\cdots&0&0\\
0&\cdots&0&0\\
\vdots&\ddots&\vdots&\vdots\\
0&\cdots&0&0
  \end{array}
\right]\\
&\triangleq&\tilde{A}+\theta_{k}\tilde{B}_{0}.
\end{eqnarray*}
Similar to Theorem 1, we can derive the following results:\\
(1) \ \ The corresponding  Riccati difference equation is
\begin{eqnarray}
H(N)\hspace{-3mm}&=&\hspace{-3mm}\tilde{H},\label{f41}\\
H_{\theta_{k-1}}(k-1,N)\hspace{-3mm}&=&\hspace{-3mm}\tilde{Q}
+\prod_{\theta_{k}}(\tilde{A}+\theta_{k}\tilde{B}_{0})'H_{\theta_{k}}(k,N)
(\tilde{A}+\theta_{k}\tilde{B}_{0})\nonumber\\
\hspace{-3mm}&&\hspace{-3mm}-\prod_{\theta_{k}}\tilde{M}'_{\theta_{k}}(k,N)\Lambda^{-1}_{\theta_{k}}(k,N)\tilde{M}_{\theta_{k}}(k,N),\label{f42}
\end{eqnarray}
in which
\begin{eqnarray}
\tilde{M}_{\theta_{k}}(k,N)\hspace{-3mm}&=&\hspace{-3mm}\tilde{B}'H_{\theta_{k}}(k,N)
\tilde{A}+\theta_{k}\bar{B}'H_{\theta_{k}}(k,N)\tilde{B}_{0},\label{f43}\\
\hspace{-3mm}&=&\hspace{-3mm}\left[
  \begin{array}{ccccc}
H^{2,1}_{k,N}A&H^{2,3}_{k,N}&\cdots&H^{2,d+1}_{k,N}&\theta_{k}H^{2,1}_{k,N}B
  \end{array}
\right],\nonumber\\
\Lambda_{\theta_{k}}(k,N)\hspace{-3mm}&=&\hspace{-3mm}R+\bar{B}'H_{\theta_{k}}(k,N)\bar{B}=R+H^{2,2}_{k,N};\label{f44}
\end{eqnarray}
(2) \ The costate is
\begin{eqnarray}
\beta_{k-1}=H_{\theta_{k-1}}(k-1,N)z(k);\label{f450}
\end{eqnarray}
(3) \ The optimal control is
\begin{eqnarray}
u^{\ast}(k)\hspace{-3mm}&=&\hspace{-3mm}-\Lambda^{-1}_{\theta_{k}}(k,N)\tilde{M}_{\theta_{k}}(k,N)z(k)\nonumber\\
\hspace{-3mm}&=&\hspace{-3mm}-(R+H^{2,2}_{k,N})^{-1}\big[
H^{2,1}_{k,N}Ax(k)+H^{2,3}_{k,N}u(k-1)+\cdots\nonumber\\
\hspace{-3mm}&&\hspace{-3mm}+H^{2,d+1}_{k,N}u(k-d+1)+\theta_{k}H^{2,1}_{k,N}Bu(k-d)
\big],\label{f45}
\end{eqnarray}
where $\tilde{H}_{\theta_{k}}(k,N)=(H^{i,j}_{k,N})_{(d+1)\times (d+1)}$ and $H^{i,j}_{k,N}$ denotes block matrix with suitable dimension.

Comparing (\ref{f14}) with (\ref{f45}), the following relationship will be derived
\begin{eqnarray}
\left\{
\begin{array}{lll}
(R+H^{2,2}_{k,N})^{-1}H^{2,1}_{k,N}A&=&
\Gamma^{-1}_{\theta_{k-1}}(N)M^{0}_{\theta_{k-1}}(N)\\
(R+H^{2,2}_{k,N})^{-1}H^{2,3}_{k,N}&=&
\Gamma^{-1}_{\theta_{k-1}}(N)M^{d}_{\theta_{k-1}}(N)\\
\ \ \ &\vdots&\\
(R+H^{2,2}_{k,N})^{-1}H^{2,d+1}_{k,N}&=&
\Gamma^{-1}_{\theta_{k-1}}(N)M^{2}_{\theta_{k-1}}(N)\\
(R+H^{2,2}_{k,N})^{-1}H^{2,1}_{k,N}\theta_{k}H^{2,1}_{k,N}B&=&
\Gamma^{-1}_{\theta_{k-1}}(N)M^{1}_{\theta_{k-1}}(N).\label{f46}
\end{array}
\right.
\end{eqnarray}
From  (\ref{f46}), the convergence of $\Gamma^{-1}_{\theta_{k-1}}(N)M^{j}_{\theta_{k-1}}(N), j=0,1,\cdots, d$ can be established from the convergence of $\tilde{H}_{\theta_{k}}(k,N)$ which can be obtained  in a similar manner with \cite{21}.

Let $\beta_{k-1}=\left[
  \begin{array}{c}
\beta^{0}_{k-1}\\
\beta^{1}_{k-1}\\
\vdots\\
\beta^{d}_{k-1}
  \end{array}
\right]$, and from (\ref{f450}) we know that
\begin{eqnarray}
\beta^{0}_{k-1}=H^{1,1}_{k-1,N}x(k)+H^{1,2}_{k-1,N}u(k-1)+\cdots+H^{1,d+1}_{k-1,N}u(k-d).\label{f47}
\end{eqnarray}

Further, from (\ref{f40}), we have that
\begin{eqnarray}
\beta^{0}_{k-1}=Qx(k)+E[A'\beta^{0}_{k}|\mathcal{F}_{k-1}],\label{f48}
\end{eqnarray}
comparing with (\ref{f5}) and (\ref{ff6}), it's easy to see that when $\beta^{0}_{N}=H$ the following relationship is satisfied, i.e.,
\begin{eqnarray}
\beta^{0}_{k-1}=\lambda_{k-1}. \label{f49}
\end{eqnarray}
Considering (\ref{f14}), (\ref{f16}), (\ref{f47}) and (\ref{f49}), by simply calculating, we can find the following relationship, i.e.,
\begin{eqnarray}
H^{1,1}_{k-1,N}-H^{1,2}_{k-1,N}\Gamma^{-1}_{\theta_{k-1}}(N)M^{0}_{\theta_{k}}(N)
=P_{\theta_{k-1}}(N).\label{f50}
\end{eqnarray}
From the above discussion, it is easy to know that $P_{\theta_{k}}(N)$  is convergent, i.e.,
\begin{eqnarray}
\lim\limits_{N\rightarrow\infty}P_{\theta_{k-1}}(N)\triangleq P_{m_{d-1}},\label{f500}
\end{eqnarray}
in which $\theta_{k-1}=m_{d-1}, k\geq d, m_{d-1}=0,1.$
In view of (\ref{f7}), we know that $(M^{0}_{\theta_{k-1}}(N))'\Gamma^{-1}_{\theta_{k-1}}(N)M^{0}_{\theta_{k-1}}(N)$ is convergent.

 Thereby, from the above discussion, it's not hard to verify that $\Gamma_{\theta_{k-1}}(N)$, $M^{i}_{\theta_{k-1}}(N),i=0,\cdots, d$,  $F^{d-j+1}_{\theta_{k-j-1}}(N)$ and $S^{j}_{\theta_{k-1}}(N),j=1,\cdots, d$ are also convergent and  (\ref{f22})-(\ref{f52}) are satisfied.

(2) \ In this part, the following inequality will be proved,
 \begin{eqnarray*}
P_{m_{d-1}}
-\sum_{s=1}^{d}(F^{d-s+1}_{m_{d-s}})'
\Gamma^{-1}_{m_{d-s-1}}
F^{d-s+1}_{m_{d-s}}>0.
\end{eqnarray*}
Following from Lemma 3 in \cite{18}, we have that there exists an integer $N_{0}$ such that $P_{\theta_{d-1}}(N_{0})
-\sum_{s=1}^{d}(F^{d-s+1}_{\theta_{d-s}}(N_{0}))'
\Gamma^{-1}_{\theta_{d-s-1}}(N_{0})
F^{d-s+1}_{\theta_{d-s}}(N_{0})> 0$.

Moreover,  the fact that
\begin{eqnarray*}
&&P_{m_{d-1}}(N)
-\sum_{s=1}^{d}(F^{d-s+1}_{m_{d-s}}(N))'
\Gamma^{-1}_{m_{d-s-1}}(N)
F^{d-s+1}_{m_{d-s}}(N)\\
&=&P_{m_{d-1}}(N-k+d)-\sum_{s=1}^{d}(F^{d-s+1}_{m_{d-s}}(N-k+d))'\\
&&\cdot\Gamma^{-1}_{m_{d-s-1}}(N-k+d)
F^{d-s+1}_{m_{d-s}}(N-k+d)
\end{eqnarray*}
is monotonically increasing with respect to $N$  yields that
\begin{eqnarray*}
&&P_{m_{d-1}}
-\sum_{s=1}^{d}(F^{d-s+1}_{m_{d-s}})'
\Gamma^{-1}_{m_{d-s-1}}
F^{d-s+1}_{m_{d-s}}\\
&=&\underset{N\rightarrow\infty}{\lim}[P_{m_{d-1}}(N)-\sum_{s=1}^{d}
(F^{d-s+1}_{m_{d-s}}(N))'\\
&&\cdot\Gamma^{-1}_{m_{d-s-1}}(N)
F^{d-s+1}_{m_{d-s}}(N)]\\
&\geq& P_{m_{d-1}}(N_{0})-\sum_{s=1}^{d}(F^{d-s+1}_{m_{d-s}}(N_{0}))'\\
&&\cdot\Gamma^{-1}_{m_{d-s-1}}(N_{0})
F^{d-s+1}_{m_{d-s}}(N_{0})\\
&>&0.
\end{eqnarray*}
The proof is completed.

\section{Proof of Theorem 3}
\emph{Proof}: ``\emph{Sufficiency}" Suppose that algebraic Riccati equation (\ref{f22})-(\ref{f27}) has a unique solution satisfying $P_{m_{d-1}}-\sum_{s=1}^{d}(F^{d-s+1}_{m_{d-s}})'
\Gamma^{-1}_{m_{d-s-1}}
F^{d-s+1}_{m_{d-s}}>0$, we will show the system (\ref{x1}) is stabilizable in the mean square sense. For this purpose, we first define Lyapunov function as
\begin{eqnarray}
\hspace{-3mm}&&\hspace{-3mm}\mathcal{L}(k)\nonumber\\
\hspace{-3mm}&=&\hspace{-3mm}E\bigg\{x'(k)P_{\theta_{k-1}}x(k)
-x'(k)\sum_{s=1}^{d}(F^{d-s+1}_{\theta_{k-s}})'
\Gamma^{-1}_{\theta_{k-s-1}}
M^{0}_{\theta_{k-s-1}}\nonumber\\
\hspace{-3mm}&&\hspace{-3mm}\cdot x(k-s)-x'(k)\sum_{s=0}^{d-1}\sum_{i=d-s}^{d}
(F^{d-i+1}_{\theta_{k-i}})'\Gamma^{-1}_{\theta_{k-i-1}}
M^{s+1-d+i}_{\theta_{k-i-1}}\nonumber\\
\hspace{-3mm}&&\hspace{-3mm}\cdot u(k-2d+s)\bigg\}.\label{f58}
\end{eqnarray}
The monotonicity and boundedness of the function $\mathcal{L}(k)$ will be illustrated.
From (\ref{f1105}), when $u(k-d)=-\Gamma^{-1}_{m_{d}}M^{0}_{m_{d}}x(k-d)
-\Gamma^{-1}_{m_{d}}\sum_{j=1}^{d}M^{j}_{m_{d}}u(k-2d+j-1)$, it yields that
\begin{eqnarray}
\hspace{-3mm}&&\hspace{-3mm}\mathcal{L}(k)-\mathcal{L}(k+1)\nonumber\\
\hspace{-3mm}&=&\hspace{-3mm}E\{x'(k)Qx(k)+u'(k-d)Ru(k-d)\}\geq0, k\geq d,\label{f59}
\end{eqnarray}
i.e., $\mathcal{L}(k)$  decreases with respect to $k$. With regard to its boundedness, considering (\ref{f020})-(\ref{f022}) given in Lemma 1, $\mathcal{L}(k)$ can be further expressed as following
\begin{eqnarray}
\mathcal{L}(k)\hspace{-3mm}&=&\hspace{-3mm}\mbox{E}\bigg\{x(k)'\Big[P_{\theta_{k-1}}
-\sum_{s=1}^{d}(F^{d-s+1}_{\theta_{k-s}})'\Gamma^{-1}_{\theta_{k-s-1}}
F^{d-s+1}_{\theta_{k-s}}\Big]x(k)\nonumber\\
\hspace{-3mm}&&\hspace{-3mm}+\sum_{s=1}^{d}\{F^{d-s+1}_{\theta_{k-s}}x(k)-E[F^{d-s+1}_{\theta_{k-s}}x(k)
|\mathcal{F}_{k-s-1}]\}'\Gamma^{-1}_{\theta_{k-s-1}}\nonumber\\
\hspace{-3mm}&&\hspace{-3mm}\cdot \{F^{d-s+1}_{\theta_{k-s}}x(k)-E[F^{d-s+1}_{\theta_{k-s}}x(k)
|\mathcal{F}_{k-s-1}]\}\bigg\}\nonumber\\
\hspace{-3mm}&\geq&\hspace{-3mm}\mbox{E}\bigg\{x(k)'\Big[P_{\theta_{k-1}}
-\sum_{s=1}^{d}(F^{d-s+1}_{\theta_{k-s}})'\Gamma^{-1}_{\theta_{k-s-1}}
F^{d-s+1}_{\theta_{k-s}}\Big]x(k)\bigg\}\nonumber\\
\hspace{-3mm}&\geq&\hspace{-3mm}0, \ \ k\geq d. \label{f60}
\end{eqnarray}
Thus, in consideration of its monotonicity, $\mathcal{L}(k)$ is convergent.

For any nonnegative integer $l$, when the both sides of (\ref{f59}) are summed up from $k=l+d$ to $k=l+N$ and letting $l\rightarrow\infty$, we can derive that
\begin{eqnarray}
\hspace{-3mm}&&\hspace{-3mm}\lim\limits_{l\rightarrow\infty}\sum_{k=l+d}^{l+N}E[x'(k)Qx(k)+u'(k-d)Ru(k-d)]\nonumber\\
\hspace{-3mm}&=&\hspace{-3mm}\lim\limits_{l\rightarrow\infty}[\mathcal{L}(l+d)-\mathcal{L}(l+N+1)]=0.\label{f61}
\end{eqnarray}
Recall that
\begin{eqnarray*}
\hspace{-3mm}&&\hspace{-3mm}\sum_{k=d}^{N}E[x'(k)Qx(k)+u'(k-d)Ru(k-d)]\\
\hspace{-3mm}&\geq &\hspace{-3mm} \mbox{E}\bigg\{x(d)'\Big[P_{\theta_{d-1}}(N)-\sum_{s=1}^{d}(F^{d-s+1}_{\theta_{d-s}}(N))'\nonumber\\
\hspace{-3mm}&&\hspace{-3mm}\cdot\Gamma^{-1}_{\theta_{d-s-1}}(N)
F^{d-s+1}_{\theta_{d-s}}(N)\Big]x(d)\bigg\}.
\end{eqnarray*}
Therefore, the following relationship can be deduced that
\begin{eqnarray*}
\hspace{-3mm}&&\hspace{-3mm}\sum_{k=l+d}^{l+N}E[x'(k)Qx(k)+u'(k-d)Ru(k-d)]\\
\hspace{-3mm}&\geq & \hspace{-3mm} \mbox{E}\bigg\{x(l+d)'\Big[P_{\theta_{l+d-1}}(l+N)-\sum_{s=1}^{d}
(F^{d-s+1}_{\theta_{l+d-s}}(l+N))'\nonumber\\
\hspace{-3mm}&&\hspace{-3mm}\cdot\Gamma^{-1}_{\theta_{l+d-s-1}}(l+N)
F^{d-s+1}_{\theta_{l+d-s}}(l+N)\Big]x(l+d)\bigg\}\\
\hspace{-3mm}&=&\hspace{-3mm}\mbox{E}\bigg\{x(l+d)'\Big[P_{\theta_{d-1}}(N)
-\sum_{s=1}^{d}(F^{d-s+1}_{\theta_{d-s}}(N))'\nonumber\\
\hspace{-3mm}&&\hspace{-3mm}\cdot\Gamma^{-1}_{\theta_{d-s-1}}(N)
F^{d-s+1}_{\theta_{d-s}}(N)\Big]x(l+d)\bigg\}\\
\hspace{-3mm}&\geq&\hspace{-3mm} 0.
\end{eqnarray*}
And from (\ref{f61}), we have that
\begin{eqnarray}
\hspace{-3mm}&&\hspace{-3mm}\lim\limits_{l\rightarrow\infty}\mbox{E}\bigg\{x(l+d)'
\Big[P_{m_{d-1}}(N)-\sum_{s=1}^{d}(F^{d-s+1}_{m_{d-s}}(N))'\nonumber\\
\hspace{-3mm}&&\hspace{-3mm}\cdot\Gamma^{-1}_{m_{d-s-1}}(N)
F^{d-s+1}_{m_{d-s}}(N)\Big]x(l+d)\bigg\}\nonumber\\
\hspace{-3mm}&=&\hspace{-3mm}0, \ \ \forall N\geq d.\label{f62}
\end{eqnarray}
In the light of the conclusion of Theorem 2, we have that there exists $N_{0}$ such that
\begin{eqnarray*}
\hspace{-4mm}&&\hspace{-4mm}P_{m_{d-1}}(N_{0})-\sum_{s=1}^{d}(F^{d-s+1}_{m_{d-s}}(N_{0}))'
\Gamma^{-1}_{m_{d-s-1}}(N_{0}))
F^{d-s+1}_{m_{d-s}}(N_{0}))\\
\hspace{-4mm}&>&\hspace{-4mm}0.
\end{eqnarray*}
Hence, (\ref{f62}) implies that
\begin{eqnarray*}
\lim\limits_{l\rightarrow\infty}\mbox{E}[x(l+d)'x(l+d)]=0.
\end{eqnarray*}
Therefore, the controller (\ref{f56}) stabilizes (\ref{x1}) in the mean-square sense.

Next we will illustrate that the cost (\ref{J2}) can be minimized by the controller (\ref{f56}).

Summing up from $k=0$ to $k=N$ on both sides of (\ref{f59}), it yields that
\begin{eqnarray}
\hspace{-3mm}&&\hspace{-3mm}E\left\{\sum_{k=0}^{N}x'(k)Qx(k)+\sum_{k=d}^{N}u'(k-d)Ru(k-d)\right\}\nonumber\\
\hspace{-3mm}&=&\hspace{-3mm}\mathcal{L}(0)-\mathcal{L}(N+1)-\sum_{k=0}^{d-1}u'(k-d)Ru(k-d)\nonumber\\
\hspace{-3mm}&&\hspace{-3mm}+\sum_{k=0}^{N}E\bigg\{\Big[u(k-d)+
\Gamma^{-1}_{m_{d}}M^{0}_{m_{d}} x(k-d)+\Gamma^{-1}_{m_{d}}\nonumber\\
\hspace{-3mm}&&\hspace{-3mm}\cdot\sum_{j=1}^{d}M^{j}_{m_{d}}u(k-2d+j-1)\Big]'
\Gamma_{m_{d}}\Big[u(k-d)+\Gamma^{-1}_{m_{d}}\nonumber\\
\hspace{-3mm}&&\hspace{-3mm}\cdot M^{0}_{m_{d}}x(k-d)+\Gamma^{-1}_{m_{d}}\sum_{j=1}^{d}M^{j}_{m_{d}}u(k-2d+j-1)\Big]\bigg\}.\label{f63}
\end{eqnarray}
Further, considering that
\begin{eqnarray*}
0\hspace{-3mm}&\leq&\hspace{-3mm}\mathcal{L}(k)\\
\hspace{-3mm}&=&\hspace{-3mm}E\bigg\{x'(k)P_{\theta_{k-1}}x(k)
-x'(k)\sum_{s=1}^{d}(F^{d-s+1}_{\theta_{k-s}})'
\Gamma^{-1}_{\theta_{k-s-1}}\nonumber\\
\hspace{-3mm}&&\hspace{-3mm}\cdot
M^{0}_{\theta_{k-s-1}}x(k-s)-x'(k)\sum_{s=0}^{d-1}\sum_{i=d-s}^{d}
(F^{d-i+1}_{\theta_{k-i}})'\Gamma^{-1}_{\theta_{k-i-1}}\nonumber\\
\hspace{-3mm}&&\hspace{-3mm}\cdot
M^{s+1-d+i}_{\theta_{k-i-1}}u(k-2d+s)\bigg\}\\
\hspace{-3mm}&\leq&\hspace{-3mm} E(x'(k)P_{\theta_{k-1}}x(k)),
\end{eqnarray*}
 and due to the fact that the system (\ref{x1}) is stabilized in the mean-square sense, therefore, $\lim\limits_{k\rightarrow\infty}E(x'(k)P_{\theta_{k-1}}x(k))=0$, i.e.,
  $\lim\limits_{k\rightarrow\infty}\mathcal{L}(k)=0$.

Let $N\rightarrow\infty$ on both sides of (\ref{f63}), then
\begin{eqnarray}
J\hspace{-3mm}&\leq&\hspace{-3mm}\mathcal{L}(0)-\sum_{k=0}^{d-1}u'(k-d)Ru(k-d)
+\sum_{k=0}^{d-1}E\bigg\{\Big[u(k-d)\nonumber\\
\hspace{-3mm}&&\hspace{-3mm}+\Gamma^{-1}_{m_{d}}M^{0}_{m_{d}} x(k-d)+\Gamma^{-1}_{m_{d}}\sum_{j=1}^{d}M^{j}_{m_{d}}u(k-2d+j-1)\Big]'\nonumber
\end{eqnarray}
\begin{eqnarray}
\hspace{-3mm}&&\hspace{-3mm}\cdot
\Gamma_{m_{d}}\Big[u(k-d)+
\Gamma^{-1}_{m_{d}}M^{0}_{m_{d}}x(k-d)+\Gamma^{-1}_{m_{d}}\sum_{j=1}^{d}M^{j}_{m_{d}}\nonumber\\
\hspace{-3mm}&&\hspace{-3mm}u(k-2d+j-1)\Big]+\sum_{k=d}^{\infty}E\bigg\{\Big[u(k-d)+
\Gamma^{-1}_{m_{d}}M^{0}_{m_{d}} x(k-d)\nonumber\\
\hspace{-3mm}&&\hspace{-3mm}+\Gamma^{-1}_{m_{d}}\sum_{j=1}^{d}M^{j}_{m_{d}}u(k-2d+j-1)\Big]'
\Gamma_{m_{d}}\Big[u(k-d)+
\Gamma^{-1}_{m_{d}}\nonumber\\
\hspace{-3mm}&&\hspace{-3mm}\cdot M^{0}_{m_{d}}x(k-d)+\Gamma^{-1}_{m_{d}}\sum_{j=1}^{d}M^{j}_{m_{d}}u(k-2d+j-1)\Big]\bigg\}.\label{f64}
\end{eqnarray}
In view of the positive definiteness of $\Gamma_{m_{d}}$, the optimal controller
to minimize (\ref{f64}) must be (\ref{f56}) and the corresponding optimal
cost is as (\ref{f57}). Therefore the proof of sufficiency is finished.

``\emph{Necessity}" From the above discussion in Theorem 2,  the following relationship is established
\begin{eqnarray*}
P_{m_{d-1}}-\sum_{s=1}^{d}(F^{d-s+1}_{m_{d-s}})'
\Gamma^{-1}_{m_{d-s-1}}
F^{d-s+1}_{m_{d-s}}>0.
\end{eqnarray*}
The uniqueness can be similarly derived from the proof in \cite{18}, so we omit it here.
% use section* for acknowledgment
%\section*{Acknowledgment}

%The authors would like to thank...

% Can use something like this to put references on a page
% by themselves when using endfloat and the captionsoff option.
\ifCLASSOPTIONcaptionsoff
  \newpage
\fi

\end{document}